\documentclass[12pt,a4paper]{article}

\usepackage{graphicx}
\usepackage{amsmath,amsfonts,amsthm}

\newtheorem{theorem}{Theorem}

\newcommand{\bfx}{\boldsymbol{x}}
\newcommand{\bn}{\boldsymbol{n}}
\newcommand{\bQ}{\boldsymbol{Q}}
\newcommand{\bq}{\boldsymbol{q}}
\newcommand{\bM}{\boldsymbol{M}}
\newcommand{\bN}{\boldsymbol{N}}

\newcommand{\bS}{\boldsymbol{S}}
\newcommand{\bsigma}{\boldsymbol{\sigma}}

\newcommand{\bW}{\boldsymbol{W}}
\newcommand{\bw}{\boldsymbol{w}}
\newcommand{\bWh}{\boldsymbol{W}_{\!h}}
\newcommand{\bWz}{\boldsymbol{W_z}}
\newcommand{\by}{{\boldsymbol{y}}}
\newcommand{\bz}{{\boldsymbol{z}}}
\newcommand{\cE}{\mathcal{E}}
\newcommand{\cmpl}{\mathrm{cmpl}}
\newcommand{\cN}{\mathcal{N}}
\newcommand{\cP}{\mathcal{P}}
\newcommand{\CR}{{\mathrm{CR}}}
\newcommand{\cT}{\mathcal{T}}
\newcommand{\diam}{\operatorname{diam}}
\newcommand{\ddiv}{\operatorname{div}}
\newcommand{\dx}[1][x]{\,\mathrm{d}#1}
\newcommand{\Hdiv}[1][\Omega]{\boldsymbol{H}(\ddiv,#1)}
\newcommand{\ICM}{\mathcal{I}_{\mathrm{CM}}}
\newcommand{\incl}{{\mathrm{LG}}}
\newcommand{\Kato}{{\mathrm{Kato}}}

\newcommand{\norm}[1]{\left\| #1 \right\|}
\newcommand{\R}{\mathbb{R}}
\newcommand{\RT}{\boldsymbol{\mathrm{RT}}}
\newcommand{\tbsigma}{\boldsymbol{\tilde\sigma}}
\newcommand{\VCR}{V^\CR}
\newcommand{\Wein}{{\mathrm{Wein}}}

\begin{document}

\title{%
Three methods for two-sided bounds of eigenvalues -- a comparison
} 
 
\author{Tom\'a\v{s} Vejchodsk\'y\\[2mm]
\parbox{\textwidth}{
\begin{center}
Institute of Mathematics, Czech Academy of Sciences\\
\v{Z}itn\'a 25, Praha 1, CZ-115\,67, Czech Republic\\
vejchod@math.cas.cz
\end{center}
}
}


\maketitle

\begin{abstract}
We compare three finite element based methods designed for two-sided bounds of eigenvalues of symmetric elliptic second order operators.
The first method is known as the Lehmann--Goerisch method.
The second method is based on Crouzeix--Raviart nonconforming finite element method.
The third one is a combination of generalized Weinstein and Kato bounds
with complementarity based estimators.
We concisely describe these methods and use them to solve three numerical examples. 
We compare their accuracy, computational performance, and generality
in both the lowest and higher order case.
\end{abstract}
 


\section{Introduction}

The standard conforming finite element discretization of symmetric second order elliptic eigenvalue problems \cite{BabOsb:1991,Boffi:2010} is very efficient and since it is a special case of the Ritz--Galerkin method, it yields natural upper bounds on the exact eigenvalues.
Interestingly, lower bounds are much more difficult to compute. The problem of lower bounds attracts attention for many decades. Lower bounds of Temple \cite{Temple1928} from 1928 and Weinstein \cite{Weinstein1937} from 1937 were generalized by Kato \cite{Kato1949} in 1949 and subsequently by Harrell \cite{Harrell1978}. In a sense optimal lower bounds are due to Lehmann \cite{Lehmann1949,Lehmann1950}. These bounds were made computational by Goerisch \cite{GoeHau1985}.
Lower bounds on eigenvalues are still a subject of active research,
see for example interesting ideas in \cite{GruOva2009,Kobayashi2015,Kuz_Rep_Guar_low_bound_smal_eig_elip_13,Repin2012}. 
The Lehmann--Goerisch method,
see \cite{BehnkeGoerish1994} and the overview in \cite{Plum1997}, is formulated in an abstract way using linear operators on Hilbert spaces and it is not straightforward to use these results in the context of the finite element method. A practical approach how to compute eigenvalue inclusions by the Lehmann--Goerisch method using the finite element method is described in \cite{BehMerPluWie2000} and we use this approach below as the first method for the comparison.

Papers \cite{LiuOis2011,LiuOis2013} propose to compute lower bounds on eigenvalues by employing conforming finite elements and computable estimates of an interpolation constant. This constant is estimated by solving a dense matrix eigenvalue problem, which may be time consuming. 
Results \cite{CarGal2014,CarGed2014} bring an idea to use nonconforming elements, where the interpolation constant is explicitly known. 
This result was subsequently improved in \cite{Liu2015} by removing the separation condition for higher eigenvalues and by providing sharper value for the interpolation constant. 
We choose the method based on nonconforming Crouzeix--Raviart elements with explicitly known interpolation constant \cite{CarGed2014,Liu2015} as the second one for the comparison.
We note that there exist several methods for lower bounds on eigenvalues based on nonconforming finite elements, see e.g. \cite{ArmDur2004,LuoLinXie:2012,AndRac:2012,HuHuaLin2014,HuHuaShe2014,LinLuoXie2014,LinXieLuoLiYan:2010,YanHanBiYu2015,YanZhaLin:2010}. However, the distinctive feature of the chosen method is that it does not require an \emph{a priori} information about the spectrum and provides guaranteed lower bounds even on rough meshes.

The third method, we present and compare, is based on a combination of 
Weinstein and Kato bounds generalized to the weak setting
with complementarity based estimators. This method was recently proposed in \cite{VejSeb2017} and it is a generalization of the approach from \cite{SebVej2014}. 

Generality of the chosen methods varies. The method based on Crouzeix--Raviart elements is the least general, 
because the needed interpolation constant is explicitly known for simple operators, such as the Laplace \cite{CarGed2014,Liu2015} and biharmonic \cite{CarGal2014} operators.
This is also the reason, why we choose the Laplace eigenvalue problem for the comparison. 
We seek eigenvalues $\lambda_i$ and eigenfunctions $u_i \neq 0$, $i=1,2,\dots$, defined in an open bounded Lipschitz domain $\Omega \subset \R^2$ such that
\begin{alignat}{2}
  \label{eq:modprob}
  - \Delta u_i &= \lambda_i u_i &\quad&\text{in }\Omega, 
  \\ \nonumber
           u_i &=0 &\quad&\text{on }\partial\Omega.
\end{alignat}
The weak formulation of this problem is based on the Sobolev space $V=H^1_0(\Omega)$ of square integrable functions with square integrable distributional derivatives and zero traces on the boundary $\partial\Omega$. 
Denoting the $L^2(\Omega)$ inner product by $(\cdot,\cdot)$,
the weak formulation reads: find $\lambda_i \in \R$ and $u_i \in V$, $u_i \neq 0$, such that
\begin{equation}
\label{eq:weakf}
  (\nabla u_i, \nabla v) = \lambda_i (u_i, v) \quad \forall v\in V.
\end{equation}
It is well known \cite{BabOsb:1991,Boffi:2010} that eigenvalues $\lambda_i$ are positive and form a countable sequence tending to infinity. We consider the natural enumeration $0 < \lambda_1 \leq \lambda_2 \leq \cdots$ and repeat these eigenvalues according to their multiplicities. The goal is to use the three chosen methods and compute lower and upper bounds for the first $m$ eigenvalues.

The subsequent Sections~\ref{se:2}--\ref{se:4} describe the three methods we compare. Section~\ref{se:sqrpipi} presents the numerical performance of the three methods on a square domain, where the analytic solution is known. Section~\ref{se:dumbbell} provides results for a dumbbell shaped domain, where the analytic solution is not available, and Section~\ref{se:ctsqr} presents the case of an asymmetric domain. Performance of higher order versions of these methods is numerically illustrated in Section~\ref{se:pdeg}. Finally, Section~\ref{se:conclusions} draws the conclusions. 

\section{The Lehmann--Goerisch method}
\label{se:2}
This section describes the Lehmann--Goerisch method \cite{GoeHau1985,Lehmann1949,Lehmann1950} and in particular its 
implementation from \cite{BehMerPluWie2000}.
The upper bound on eigenvalues is obtained by the standard conforming finite element method. 
For simplicity, the domain $\Omega$ is considered to be a polygon. The standard finite element triangulation of $\Omega$ is denoted by $\cT_h$ and the finite element space of order $k$ is defined as
\begin{equation}
 \label{eq:defVh}
  V_h = \{ v_h \in V : v_h|_K \in P_k(K) \quad \forall K \in \cT_h \},
\end{equation}
where $P_k(K)$ is the space of polynomials of degree at most $k$ on the triangle $K$.
The finite element approximation of problem \eqref{eq:modprob} consists of seeking eigenvalues $\Lambda_{h,i} \in \R$ and corresponding eigenfunctions $u_{h,i} \in V_h$ such that
\begin{equation}
 \label{eq:fem}
  ( \nabla u_{h,i}, \nabla v_h) = \Lambda_{h,i} (u_{h,i}, v_h) \quad \forall v_h \in V_h.
\end{equation}
Approximate eigenvalues are naturally ordered $0 < \Lambda_{h,1} \leq \Lambda_{h,2} \leq \cdots$ and repeated according to their multiplicities.
It is well known \cite{BabOsb:1991,Boffi:2010} that 
for sufficiently smooth eigenfunctions
the convergence order of approximate eigenvalues $\Lambda_{h,i}$ is $2k$ as the mesh size $h$ tends to zero 
and that they provide upper bounds on the exact eigenvalues:
$\lambda_i \leq \Lambda_{h,i}$ for all $i=1,2,\dots, \operatorname{dim} V_h$.

Lower bounds on eigenvalues are based on the result provided in \cite[Theorem~2.1]{BehMerPluWie2000}. 
For the readers' convenience, we present this result here as Theorem~\ref{th:BMPW}. Note that $\bW = \Hdiv$ denotes the space of square integrable vector fields with square integrable divergence. 
\begin{theorem}\label{th:BMPW}
Let $(\tilde u_i, \tbsigma_i) \in V \times \bW$, $i=1,2,\dots,n$, and $\rho >0$, $\gamma > 0$ be arbitrary.
Define matrices $\bM, \bN \in \R^{n \times n}$ with entries
\begin{align*}
  \bM_{ij} &= (\nabla \tilde u_i, \nabla \tilde u_j) + (\gamma - \rho) (\tilde u_i, \tilde u_j), \\
  \bN_{ij} &= (\nabla \tilde u_i, \nabla \tilde u_j) + (\gamma - 2\rho) (\tilde u_i, \tilde u_j)
           + \rho^2 (\tbsigma_i,\tbsigma_j) 
           + (\rho^2/\gamma) (\tilde u_i + \ddiv\tbsigma_i,\tilde u_j + \ddiv\tbsigma_j).
\end{align*}
Suppose, that the matrix $\bN$ is positive definite, and let
$$
  \mu_1 \leq \mu_2 \leq \dots \leq \mu_n
$$
be eigenvalues of the generalized eigenvalue problem 
\begin{equation}
  \label{eq:MNproblem}
  \bM \by_i = \mu_i \bN \by_i, \quad i=1,2,\dots,n.
\end{equation}
Then, for all $i$ such that $\mu_i < 0$, the interval 
$$
  [ \rho - \gamma - \rho/(1-\mu_i), \rho - \gamma)
$$
contains at least $i$ eigenvalues of the continuous problem \eqref{eq:weakf}.
\end{theorem}

Functions $\tilde u_i \in V$ and $\tbsigma_i \in \bW$ in Theorem~\ref{th:BMPW} are in general arbitrary, but in order to obtain accurate bounds, they should approximate the exact eigenfunction $u_i$ and the corresponding flux $(\lambda_i+\gamma)^{-1} \nabla u_i$, respectively. The natural choice for $\tilde u_i$ is the finite element approximation $u_{h,i}$. 
We will follow \cite{BehMerPluWie2000} and choose the flux $\tbsigma_i$ as
a solution $\bsigma_{h,i}$ of a saddle point problem solved by the mixed finite element method. Note that this idea is closely related to approaches called dual finite elements in \cite{HasHla:1976,Hla:1978,HlaKri:1984} and the complementarity technique in \cite{Complement:2010,systemaee:2010}.
 Denoting by 
\begin{equation}
  \label{eq:RTk}
  \RT_k(K) = [P_k(K)]^2 \oplus \bfx P_k(K)
\end{equation}
the local Raviart--Thomas space on the element $K \in \cT_h$,
the flux reconstruction $\bsigma_{h,i}$ is sought in the global space 
$$
  \bWh = \{ \bsigma_h \in \Hdiv : \bsigma_h|_K \in \RT_k(K) \quad\forall K \in \cT_h \},
$$
and the Lagrange multipliers in
$$
  Q_h = \{ q_h \in L^2(\Omega) : q_h |_K \in P_k(K) \quad\forall K \in \cT_h \}.
$$
The mixed finite element problem then reads: find $(\bsigma_{h,i},q_{h,i}) \in \bWh \times Q_h$ such that
\begin{alignat}{2}
  \label{eq:sig1}
  (\bsigma_{h,i}, \bw_h) + (q_{h,i}, \ddiv \bw_h) &= 0 &\quad&\forall \bw_h \in \bWh, 
  \\
  \label{eq:sig2}
  (\ddiv \bsigma_{h,i}, \varphi_h ) &= (- u_{h,i}, \varphi_h) &\quad &\forall \varphi_h \in Q_h,
\end{alignat}
where $u_{h,i} \in V_h$ is the finite element approximation \eqref{eq:fem}.

\begin{figure}
\centerline{
\includegraphics[height=0.27\textwidth]{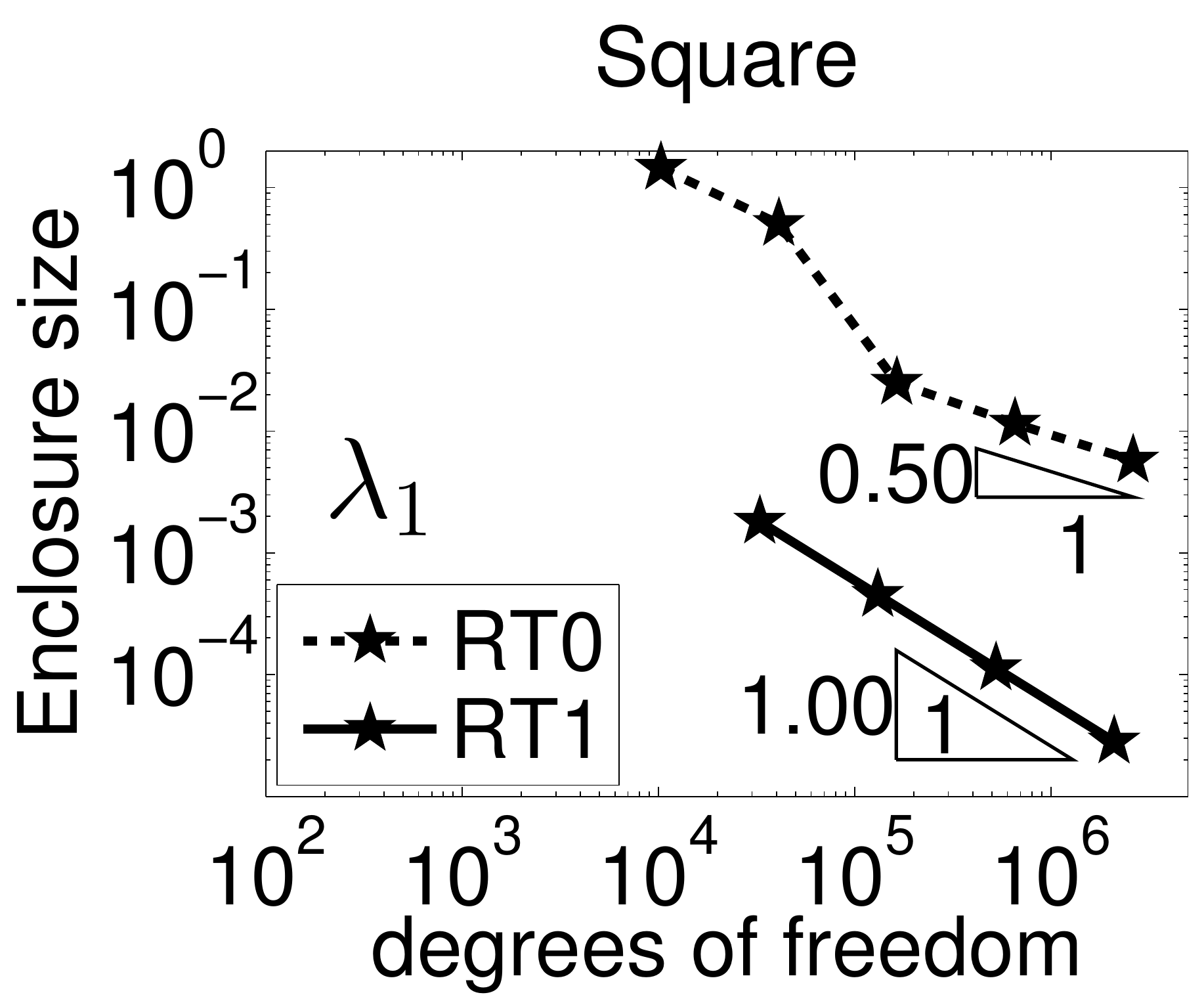}\ %
\includegraphics[height=0.27\textwidth]{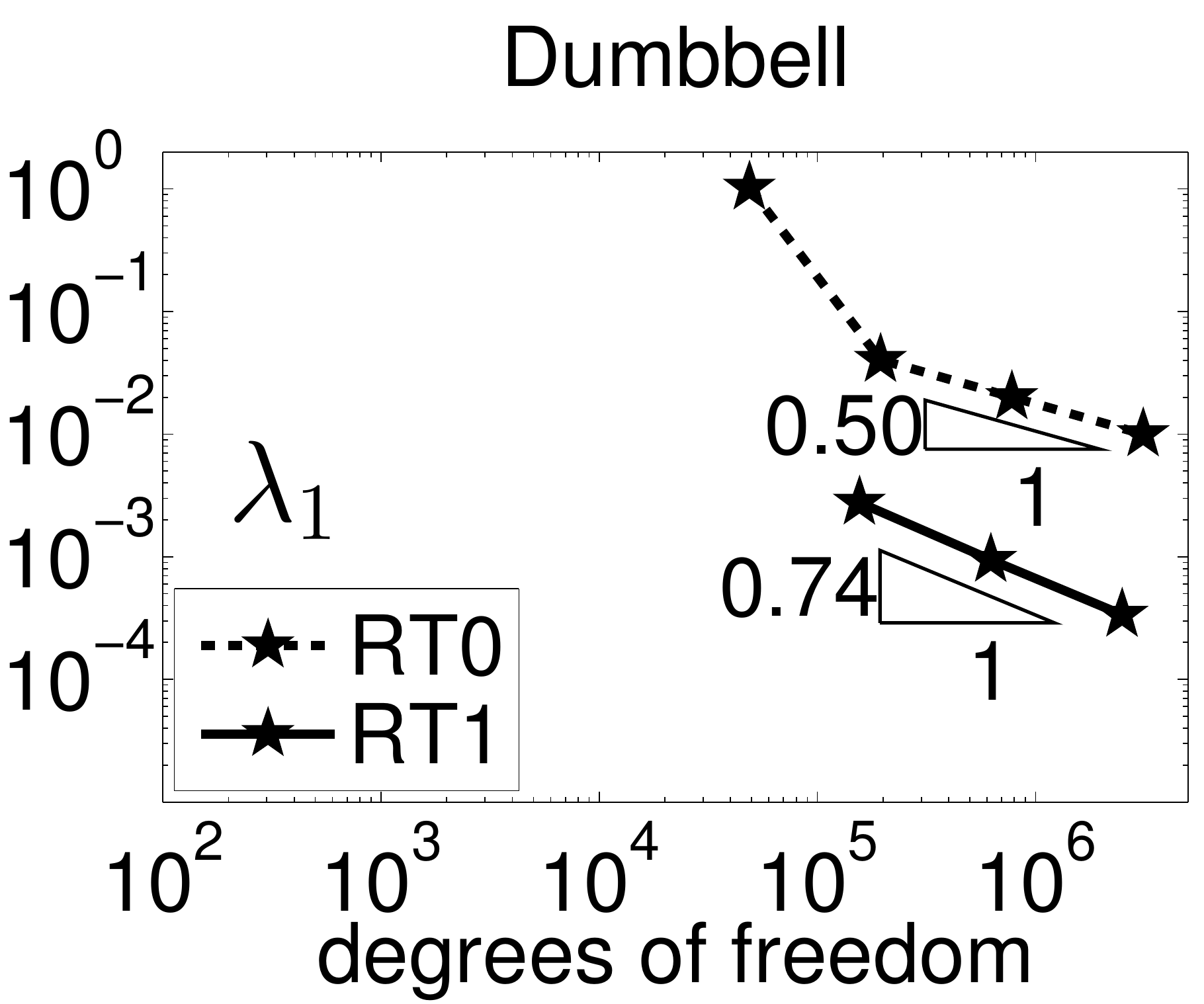}\ %
\includegraphics[height=0.27\textwidth]{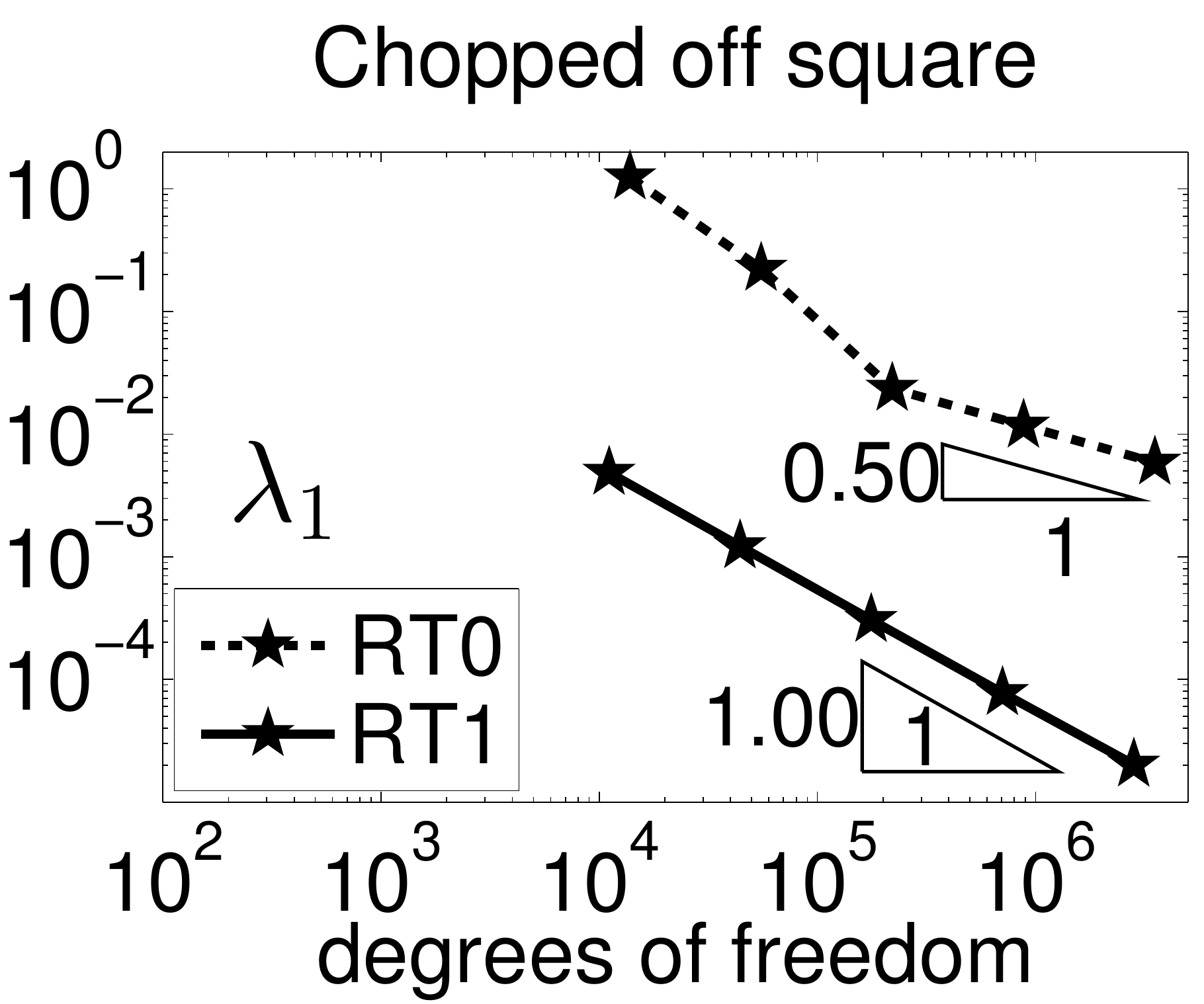}%
}
\caption{\label{fi:RT0RT1}
Comparison of choices $\RT_0(K)$ and $\RT_1(K)$ in the definition of $\bWh$ for the lowest order case ($k=1$). Figures show convergence curves of the Lehmann--Goerisch method for the first eigenvalue in the square (left), dumbbell shaped domain (middle), and chopped off square (right). For more information about these examples see Sections~\ref{se:sqrpipi}--\ref{se:ctsqr} below.
}
\end{figure}

Paper \cite{BehMerPluWie2000} proposes to choose $\RT_{k-1}(K)$ in definitions of spaces $\bWh$
(and consequently $P_{k-1}(K)$ in the definition of $Q_h$) and \cite[Remark 2.5]{BehMerPluWie2000} justifies $O(h^k)$ convergence of the resulting lower bound.
However, in test examples the optimal $O(h^{2k})$ convergence rate for regular eigenfunctions is observed for the choice 
$\RT_k(K)$, see also \cite{AinVej2014}.
This is illustrated in Figure~\ref{fi:RT0RT1} for the lowest order case ($k=1$), where the $O(h)$ rate corresponds to slope $0.5$ and $O(h^2)$ to slope $1$. 
%

To obtain lower bounds on eigenvalues based on Theorem~\ref{th:BMPW}, an \emph{a~priori} information about eigenvalues is needed. Namely,
if $\rho-\gamma \leq \lambda_{s+1}$ for some index $s \geq 1$ then Theorem~\ref{th:BMPW} provides lower bounds
\begin{equation}
\label{eq:llowinc}
\rho - \gamma - \rho/(1-\mu_i) \leq \lambda_{s+1 - i} \quad \forall i=1,2,\dots,\min\{s,n\}.
\end{equation}
Thus, if an \emph{a priori} lower bound on at least one exact eigenvalue is known then lower bounds on eigenvalues below this one can be computed by \eqref{eq:llowinc}.

Numerical examples below use an \emph{a priori} known lower bounds $\underline{\lambda}_{m+1} = \ell^\CR_{m+1}$ on $\lambda_{m+1}$ computed by the method based on Crouzeix--Raviart elements\footnote{Let me thank the anonymous referee for pointing this idea out.} described in Section~\ref{se:3}. Utilizing this information, accurate lower bounds on the first $m$ eigenvalues are obtained as follows.
\begin{enumerate}
\item 
Compute the standard finite element approximations \eqref{eq:fem} of the first $m$ eigenpairs $(\Lambda_{h,i}, u_{h,i}) \in \R \times V_h$, $i=1,2,\dots,m$. This provides upper bounds $\Lambda_{h,i}$, $i=1,2,\dots,m$, on the exact eigenvalues.
\item
Find $\bsigma_{h,i} \in \bWh$ for $i=1,2,\dots,m$ by solving \eqref{eq:sig1}--\eqref{eq:sig2}. 
\item
Apply Theorem~\ref{th:BMPW} with $\tilde u_i = u_{h,i}$, $\tbsigma_i = \bsigma_{h,i}$, $i=1,2,\dots,m$,
$\gamma = \norm{u_{h,m} + \ddiv \bsigma_{h,m}}_{L^2(\Omega)}$ and 
$\rho = \underline{\lambda}_{m+1} + \gamma$.
Assemble matrices $\bM$ and $\bN$ and find eigenvalues $\mu_1 \leq \mu_2 \leq \cdots \leq \mu_m$.
Check assumptions of the positive definiteness of $\bN$ and of the negativity of $\mu_i$.
If they are satisfied then compute lower bound \eqref{eq:llowinc}. In particular, choose $s=m$ and $i=m+1-j$ in \eqref{eq:llowinc} and get
$$
  \ell^\incl_j = \rho - \gamma - \rho/\left(1-\mu_{m+1-j}\right) \leq \lambda_j,
  \quad j = 1,2,\dots,n.
$$
%
%
\end{enumerate}

Note that the description of the Lehmann--Goerisch method provided here is tailored to the test problem \eqref{eq:modprob}. In \cite{BehMerPluWie2000}, the method is described for the Laplace eigenvalue problem with mixed homogeneous Dirichlet and Neumann boundary conditions. However, it is clearly not limited to such simple problems and can be straightforwardly generalized to problems with reaction terms and with variable diffusion and reaction coefficients.

\section{The method based on Crouzeix--Raviart elements}
\label{se:3}
This section describes the method from \cite{CarGed2014,Liu2015}. We will refer to it as the CR method, because it is based on Crouzeix--Raviart finite elements. The lowest-order variant is introduced, because a higher order version is not available. Considering the triangulation $\cT_h$ of $\Omega$ as above and the space of piecewise affine and in general discontinuous functions
$$
  \cP_1(\cT_h) = \{ v_h \in L^2(\Omega) : v_h|_K \in P_1(K) \},
$$ 
we define the standard Crouzeix--Raviart finite element space $\VCR_h \subset \cP_1(\cT_h)$ of functions $v_h \in \cP_1(\cT_h)$ that are continuous in midpoints of all interior edge of $\cT_h$ and vanish at midpoints of all boundary edges of $\cT_h$.


The Crouzeix--Raviart approximate eigenpairs $(\lambda_{h,i}^\CR,u_{h,i}^\CR) \in \R \times \VCR_h$, $u_{h,i}^\CR \neq 0$,
of problem \eqref{eq:modprob} are defined by the relation
\begin{equation}\label{eq:CRprob}
  (\nabla u_{h,i}^\CR, \nabla v_h) = \lambda_{h,i}^\CR (u_{h,i}^\CR, v_h) \quad \forall v_h \in \VCR_h.
\end{equation}
Approximate eigenvalues $\lambda_{h,i}^\CR$ are proved to be below the exact eigenvalues $\lambda_i$ for sufficiently fine meshes, but for very rough meshes it is not difficult to construct an example such that $\lambda_{h,i}^\CR$ is above $\lambda_i$, see \cite{CarGed2014}. The idea of explicit estimates of the interpolation constant enables us to construct simple lower bounds on exact eigenvalues for arbitrary meshes. 
It is proved in \cite[Theorem 3.2 and 5.1]{CarGed2014} and \cite[Theorem 2.1]{Liu2015} that 
\begin{equation}
  \ell^\CR_i \leq \lambda_i 
\quad\text{for}\quad
\label{eq:llowCR}
  \ell^\CR_i = \frac{\lambda_{h,i}^\CR}{1 + \kappa^2 \lambda_{h,i}^\CR h_{\mathrm{max}}^2},
  \quad \forall i=1,2,\dots,
\end{equation}
where 
$h_{\mathrm{max}} = \max_{K \in \cT_h} \diam K$ is the largest of all diameters of elements in the triangulation $\cT_h$ and $\kappa^2$ is a universal constant.
We use the optimal bound $\kappa \leq 0.1893$ from \cite{Liu2015} in the subsequent numerical examples.

Concerning the upper bound on eigenvalues, we can well use the standard conforming finite element approximations given by \eqref{eq:fem}. This would, however, mean to solve one more matrix eigenvalue problem.
Instead, we can use the well known idea \cite{LuoLinXie:2012,CarGed2014} of a conforming quasi-interpolation of the already computed nonconforming eigenfunction $\tilde u_{h,i}^\CR$.
The common Oswald quasi-interpolation operator \cite{Oswald1994} would be a straightforward approach to use, but it does not provide as accurate results as the quasi-interpolation proposed in \cite{CarMer2013} and used in \cite{CarGed2014}. 
This quasi-interpolation operator is defined as $\ICM : \VCR_h \to V_h^*$,
where $V_h^* = \{ v_h \in V : v_h|_K \in P_1(K) \ \forall K \in \cT_h^* \}$ and $\cT_h^*$ is the uniform (red) refinement of the triangulation $\cT_h$ such that all triangles in $\cT_h$ are refined into four similar subtriangles of $\cT_h^*$.
If $\cN_h$ stands for the set of all vertices of the triangulation $\cT_h$ and $\cE_h$ for the set of all edges in $\cT_h$ and if $v^\CR_h \in \VCR_h$ is arbitrary then
$$
\left(\ICM v^\CR_h \right) (\bz) = \left\{ \begin{array}{ll}
0 & \text{if } \bz \text{ lies on } \partial\Omega, \\
v^\CR_h(\bz) & \text{if } \bz \text{ is the midpoint of an edge } \gamma \in \cE_h,\ \gamma\not\subset\partial\Omega, \\
v_{\textrm{min}}(\bz) & \text{if } \bz \in \cN_h\setminus\partial\Omega.
\end{array} \right.
$$
The function $v_{\textrm{min}}$ is determined by a one-dimensional minimization on the patch $\omega_\bz^*$ of elements from $\cT_h^*$ sharing the vertex $\bz$. We set
$V_\bz = \{ v_h \in C(\overline{\omega}_\bz^*) : v_h|_K \in P_1(K) \text{ for all } K\in\cT_h^*,\ K\subset\overline{\omega}_\bz^*, \text{ and } v_h = v^\CR_h \text{ on } \partial \omega_\bz^* \}$ and determine $v_{\textrm{min}} \in V_\bz$ as the unique minimizer of
$$
  \min_{v_h \in V_\bz} \norm{\nabla v^\CR_h - \nabla v_h }_{L^2(\omega_\bz^*)}.
$$
Then Rayleigh quotients constructed from $u_{h,i}^* = \ICM \tilde u_{h,i}^\CR$ provide upper bounds in a standard way. 

In particular, since the goal is to obtain lower and upper bounds for the first $m$ eigenvalues, we proceed according to the following algorithm.
\begin{enumerate}
\item Solve the Crouzeix--Raviart eigenvalue problem \eqref{eq:CRprob} for $i=1,2,\dots,m$.
\item Use these approximations and compute lower bounds $\ell^\CR_i$ by \eqref{eq:llowCR} for $i=1,2,\dots,m$.
\item Construct the uniform (red) refinement $\cT_h^*$ of the mesh $\cT_h$ and interpolants $u_{h,i}^* = \ICM \tilde u_{h,i}^\CR$ for $i=1,2,\dots,m$.
\item Using the standard Ritz--Galerkin method, assemble matrices $\bS, \bQ \in \R^{m\times m}$ with entries $\bS_{j,k} = (\nabla u_{h,j}^*, \nabla u_{h,k}^*)$ and $\bQ_{j,k} = (u_{h,j}^*, u_{h,k}^*)$
and solve the matrix eigenvalue problem 
$$
  \bS \by_i = \Lambda^*_i \bQ \by_i, \quad i=1,2,\dots,m,
$$
see also \cite{LuoLinXie:2012}. Sort these eigenvalues such that $\Lambda^*_1 \leq \Lambda^*_2 \leq \cdots \leq \Lambda^*_m$ and obtain upper bounds
$$
  \lambda_i \leq \Lambda^*_i \quad \text{for } i=1,2,\dots,m.
$$
\end{enumerate}

\section{The complementarity method}
\label{se:4}
This section describes the method introduced in \cite{VejSeb2017} which is refer to as the complementarity method. It is based on the standard conforming finite element approximation \eqref{eq:fem}. Lower bounds on eigenvalues are obtained by 
a combination of generalized Weinstein \cite{Weinstein1937} and Kato \cite{Kato1949} bounds.
A crucial point is that the energy norm of the weak representative of the residual is bounded 
using the complementarity technique and local flux reconstruction \cite{BraSch:2008}, see also \cite{DolErnVoh2016,ErnVoh2013}. 

The triangulation $\cT_h$ of the domain $\Omega$ is considered as above and the finite element approximate eigenpair $(\Lambda_{h,i}, u_{h,i}) \in \R \times V_h$ as in \eqref{eq:fem}.
Based on the gradient $\nabla u_{h,i}$ of the approximate eigenvector, a suitable flux $\bq_{h,i} \in \Hdiv$ is constructed. This flux is computed by solving small mixed finite element problems on patches of elements sharing a single vertex. Let $\bz \in \cN_h$ be a vertex in $\cT_h$ and let $\cT_\bz$ be the set of those elements in $\cT_h$ that $\bz$ is one of their vertices. The patch of elements sharing the vertex $\bz$ is denoted by $\omega_\bz = \operatorname{int} \bigcup \{ K : K\in \cT_\bz\}$. 
If $\bz$ is an interior vertex then we set $\Gamma_{\omega_\bz}^\mathrm{ext} = \partial\omega_\bz$
and if $\bz$ lies on the boundary $\partial\Omega$ then we set $\Gamma_{\omega_\bz}^\mathrm{ext} = \partial\omega_\bz\setminus\partial\Omega$.
Denoting by $\bn_{\omega_\bz}$ the unit outward facing normal vector to $\partial\omega_\bz$,
spaces 
$$
  \bWz = \left\{ \bw_h \in \Hdiv[\omega_\bz] : \bw_h|_K \in \RT_k(K)\ \forall K \in \cT_\bz 
  \quad\text{and}\quad
  \bw_h \cdot \bn_{\omega_\bz} = 0 \text{ on } \Gamma_{\omega_\bz}^\mathrm{ext} \right\}
$$
and 
$$
  \cP_k^*(\cT_\bz) = \left\{ \begin{array}{ll}
    \{ v_h \in \cP_k(\cT_\bz): \int_{\omega_\bz} v_h \dx = 0 \} \quad \text{for interior vertices } \bz \in \cN_h\setminus\partial\Omega, \\
    \cP_k(\cT_\bz) \quad \text{for boundary vertices } \bz \in \cN_h \cap \partial\Omega.
    \end{array} \right.
$$
are defined.
Recall that spaces $\RT_k(K)$ were introduced in \eqref{eq:RTk} and $\cP_k(\cT_\bz)$ stands for the space of piecewise polynomials of degree at most $k$ that are in general discontinuous.
Further, $\psi_\bz$ stands for the standard piecewise affine and continuous finite element hat function associated with the vertex $\bz \in \cN_h$. Function $\psi_\bz$ has value one at $\bz$ and vanishes at all other vertices of the triangulation $\cT_h$. 
We also introduce the residual
$$
  r_{\bz,i} = \Lambda_{h,i} \psi_\bz u_{h,i} - \nabla \psi_\bz \cdot \nabla u_{h,i}.
$$

The flux reconstruction $\bq_{h,i} \in \Hdiv$ is then defined as
\begin{equation}
  \label{eq:qh}
  \bq_{h,i} = \sum_{\bz \in \cN_h} \bq_{\bz,i},
\end{equation}
where $\bq_{\bz,i} \in \bWz$ together with $d_{\bz,i} \in \cP_k^*(\cT_\bz)$ solves the mixed finite element problem
\begin{alignat}{2}
  \label{eq:locprob1}
  (\bq_{\bz,i}, \bw_h)_{\omega_\bz} - (d_{\bz,i}, \ddiv \bw_h)_{\omega_\bz} &= (\psi_\bz \nabla u_{h,i}, \bw_h)_{\omega_\bz}
  &\quad &\forall \bw_h \in \bWz, 
  \\
  \label{eq:locprob2}
  -(\ddiv \bq_{\bz,i}, \varphi_h)_{\omega_\bz} &= (r_{\bz,i}, \varphi_h)_{\omega_\bz} 
  &\quad &\forall \varphi_h \in \cP_k^*(\cT_\bz).
\end{alignat}
The reconstructed flux $\bq_{h,i}$ is used to define the error estimator
\begin{equation}
  \label{eq:eta}
  \eta_i = \norm{\nabla u_{h,i} - \bq_{h,i}}_{L^2(\Omega)}.
\end{equation}
This is finally used to define two lower bounds on eigenvalues
\begin{align}
\label{eq:Wein}
  \ell^\Wein_i &= \frac{1}{4} \left( -\eta_i + \sqrt{\eta_i^2 + 4\Lambda_{h,i}} \right)^2,
\\
  \label{eq:Kato}
  \ell^\Kato_{(s),i} &= \Lambda_{h,i} \left( 1 + \nu \lambda_{h,i} \sum_{j=i}^s \frac{\eta_j^2}{\Lambda_{h,j}^2 (\nu - \Lambda_{h,j} ) } \right)^{-1},
\end{align}
where $\nu$ is assumed to be greater than $\Lambda_{h,s}$ for some index $s \geq i$. 
The generalized Weinstein bound $\ell^\Wein_i$ is proved \cite[Theorem~2.3]{VejSeb2017} to satisfy
$$
  \ell^\Wein_i \leq \lambda_i
$$
provided the closeness assumption $\sqrt{\lambda_{i-1}\lambda_i} \leq \Lambda_{h,i} \leq \sqrt{\lambda_i \lambda_{i+1}}$ holds.
Similarly, by \cite[Theorem~2.5]{VejSeb2017} the generalized Kato bound $\ell^\Kato_{(s),i}$ satisfies
$$
  \ell^\Kato_{(s),i} \leq \lambda_i 
$$
if $\lambda_{s-1} \leq \Lambda_{h,s} < \nu \leq \lambda_{s+1}$.
The first inequality in this chain is automatically satisfied, because the finite element eigenvalue $\Lambda_{h,s}$ is an upper bound to $\lambda_s$. The second inequality can easily be verified, because both $\Lambda_{h,s}$ and $\nu$ are available. Thus, it is only the assumption $\nu \leq \lambda_{s+1}$ that represents the needed \emph{a priori} information about the spectrum. Namely, $\nu$ is an \emph{a priori} known lower bound on an exact eigenvalue in the exactly same way as $\rho-\gamma$ is the \emph{a priori} known lower bound needed in the Lehmann--Goerisch method, see the assumption above \eqref{eq:llowinc}.

Lower bound \eqref{eq:Wein} has a suboptimal rate of convergence, therefore the bound \eqref{eq:Kato} is preferred.
On the other hand, the bound \eqref{eq:Wein} is quite robust. In all numerical experiments we performed it provided lower bounds on eigenvalues. Even in cases where the relative closeness assumption was clearly not satisfied.
Therefore, we propose to combine bounds \eqref{eq:Wein} and \eqref{eq:Kato} in a recursive manner 
and compute lower bounds on the first $m$ eigenvalues using an \emph{a priori} known lower bound $\underline{\lambda}_{m+1}$ on $\lambda_{m+1}$. 
As in the Lehmann--Goerisch method, we choose $\underline{\lambda}_{m+1} = \ell^\CR_{m+1}$ computed by the CR method. The detailed algorithm reads as follows.
\begin{enumerate}
\item
Compute the standard finite element approximations $(\Lambda_{h,i}, u_{h,i}) \in \R \times V_h$, $i=1,2,\dots,m$, of the first $m$ eigenpairs according to \eqref{eq:fem}. This provides upper bounds $\Lambda_{h,i}$, $i=1,2,\dots,m$, on the exact eigenvalues.
\item 
For all $i=1,2,\dots,m$ and all vertices $\bz \in\cN_h$, solve local patch problems \eqref{eq:locprob1}--\eqref{eq:locprob2},
construct the flux $\bq_{h,i}$ as in \eqref{eq:qh}, and compute $\eta_i$ by \eqref{eq:eta}.
\item
Evaluate lower bounds $\ell^\Wein_i$ using \eqref{eq:Wein} for $i=1,2,\dots,m$
and set 
$\ell^\Kato_{m+1} = \underline{\lambda}_{m+1}$.
\item
For $n=m,m-1,\dots,1$ do the following steps.
\begin{itemize}
  \item Set $\nu = \ell^\Kato_{n+1}$. 
  \item If $\Lambda_{h,n} < \nu$ compute $\ell^\Kato_{(n),i}$ by \eqref{eq:Kato} for $i=1,2,\dots,n$.
  \item The lower bound on $\lambda_n$ is
  $$
    \ell^\Kato_n = \max \{ \ell^\Kato_{(n),n}, \ell^\Kato_{(n+1),n}, \dots, \ell^\Kato_{(m),n} \},
  $$
  where lower bounds $\ell^\Kato_{(j),n}$ for $j=n,n+1,\dots,m$ are considered only if they are defined.
\end{itemize}
\item As the final lower bound on $\lambda_i$ set
  \begin{equation}
    \label{eq:compl}
    \ell^\cmpl_i = \max \{ \ell^\Wein_i, \ell^\Kato_i \},
    \quad i = 1,2,\dots,m.
  \end{equation}
\end{enumerate}

The less accurate but sometimes easily available Weinstein bound $\underline{\lambda}_{m+1} = \ell^\Wein_{m+1}$ can also be used as the needed \emph{a priori} lower bound for both Kato and Lehmann--Goerisch methods. In this case, numerical experiments show (results not presented) that almost the same accuracy as for the choice $\underline{\lambda}_{m+1} = \ell^\CR_{m+1}$ is reached on sufficiently fine meshes, but the convergence is delayed. 

\section{Numerical results -- a square domain}
\label{se:sqrpipi}
This section computes two-sided bounds of the first $m=10$ eigenvalues
of problem \eqref{eq:modprob} in a square $\Omega = (0,\pi)^2$.
These bounds are computed by using the three methods described in Sections~\ref{se:2}--\ref{se:4} and their numerical performance is compared.
In order to obtain results comparable with the CR method, the order $k=1$ is chosen for both the Lehmann--Goerisch and complementarity methods.
The exact eigenvalues and eigenfunctions are in this case well known to be 
$$
  \lambda_{i,j} = i^2 + j^2, \quad u_{i,j}(x,y) = \sin(ix) \sin(jy), \quad i,j=1,2,\dots,
$$
see the fifth column of Table~\ref{ta:sq_eig}.
Notice that four out of the first six distinct eigenvalues are doubled.

\begin{table}
\centerline{
\begin{tabular}{cccccc}
 & lower bound       & lower bound & lower bound     & exact      & upper bound    \\ 
 & Lehmann--Goerisch & CR method   & complementarity & value & Ritz--Galerkin \\ \hline
$\lambda_{1}$ &   1.99999574 &   1.99999042 &   1.99999791 &  2  &  2.00000266 \\ 
$\lambda_{2}$ &   4.99994696 &   4.99994719 &   4.99992537 &  5  &  5.00001360 \\ 
$\lambda_{3}$ &   4.99994696 &   4.99994719 &   4.99996801 &  5  &  5.00001720 \\ 
$\lambda_{4}$ &   7.99967683 &   7.99984672 &   7.99982872 &  8  &  8.00003922 \\ 
$\lambda_{5}$ &   9.99918816 &   9.99981070 &   9.99951065 & 10  &  10.0000555 \\ 
$\lambda_{6}$ &   9.99918824 &   9.99981070 &   9.99977495 & 10  &  10.0000614 \\ 
$\lambda_{7}$ &   12.9970717 &   12.9996149 &   12.9991638 & 13  &  13.0000903 \\ 
$\lambda_{8}$ &   12.9970717 &   12.9996149 &   12.9995981 & 13  &  13.0001125 \\ 
$\lambda_{9}$ &   16.9688191 &   16.9994843 &   16.9937086 & 17  &  17.0001645 \\ 
$\lambda_{10}$&   16.9688191 &   16.9994843 &   16.9969011 & 17  &  17.0001686 \\ 
\end{tabular}
}
\caption{\label{ta:sq_eig}
Square domain. Eigenvalue bounds computed on the finest meshes. The mesh sizes of the finest meshes were $h = 0.0123$ for the Lehmann--Goerisch method, $h=0.0061$ for the CR method, and $h=0.0031$ for the complementarity method.
}
\end{table}

\begin{figure}
\centerline{
\includegraphics[height=32mm]{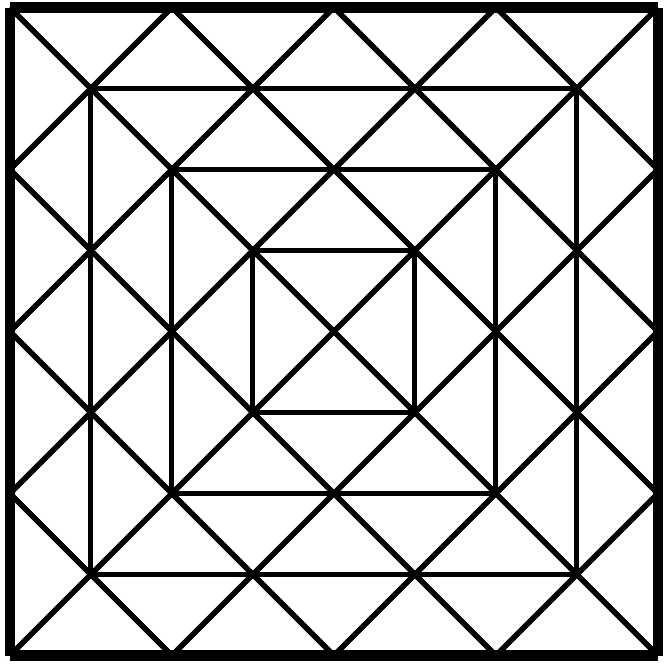}\qquad
\includegraphics[height=32mm]{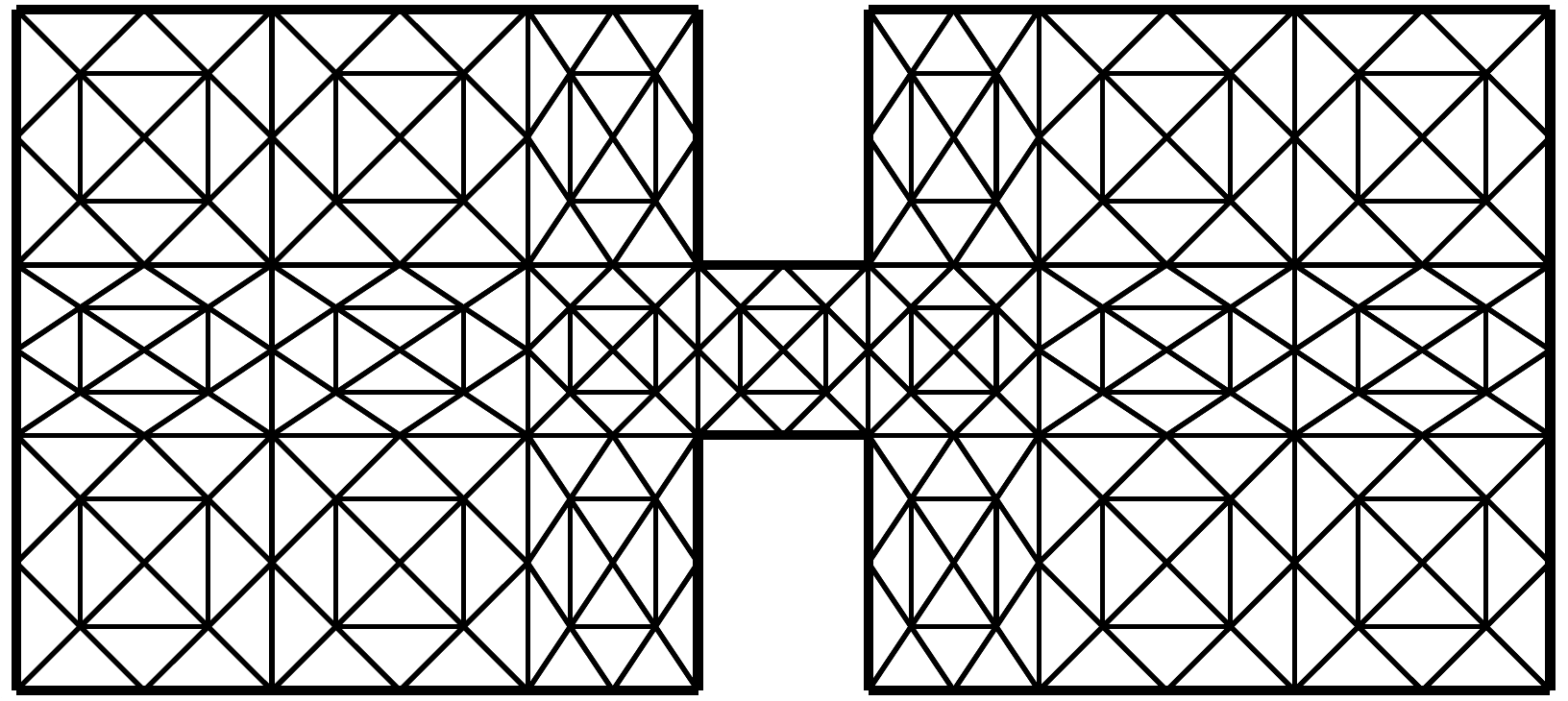}\qquad
\includegraphics[height=32mm]{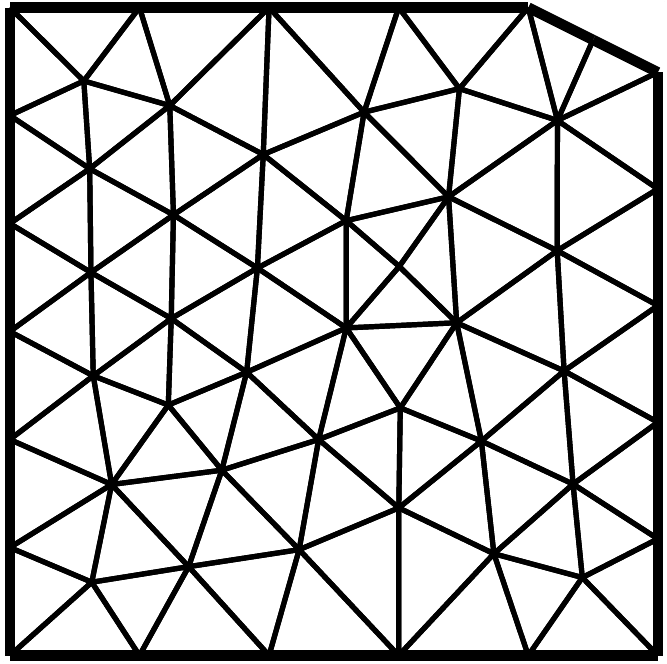}
}
\caption{\label{fi:sq_initmesh}
Initial triangulations of the square domain (left), the dumbbell shaped domain (middle), and the chopped off square (right).
}
\end{figure}

For all three methods the same triangulations are used.
The initial rough triangulation is shown in Figure~\ref{fi:sq_initmesh} (left).
Then a sequence of nested meshes is obtained by successive uniform (red) refinement. 
This means that each triangle of the original mesh is refined into four similar subtriangles.

Given a mesh $\cT_h$ of this sequence, the lower and upper bounds of eigenvalues are computed by all three methods. 
In particular, the required \emph{a~priori} lower bounds $\underline{\lambda}_{m+1}$ on $\lambda_{m+1}$ is computed by the CR method on the current mesh $\cT_h$ for both the Lehmann--Goerisch and complementarity methods.
The three methods, however, use different types of finite elements and consequently, they require to solve matrix problems of different sizes. Therefore, we compare their accuracy with respect to the number of degrees of freedom corresponding to the largest matrix problem that has to be solved.
To be more specific, we choose $N^\mathrm{mix} = \dim \bWh + \dim Q_h$ to be the reference number of degrees of freedom for the Lehmann--Goerisch method, because the mixed finite element problem \eqref{eq:sig1}--\eqref{eq:sig2} is considerably larger than problem \eqref{eq:fem} and its solution consumes most of the computational time.
For the CR method, we naturally choose $N^\CR = \dim \VCR_h$ as the reference number of degrees of freedom, although the upper bound computed by this method requires the refined mesh $\cT_h^*$ and the interpolation $\ICM \tilde u_{h,i}^\CR$ with more degrees of freedom then $N^\CR$. This interpolation, however, can be done locally on patches and does not consume the majority of the computational time. Finally, for the complementarity method we choose $N^\mathrm{conf} = \dim V_h$ as the reference number of degrees of freedom. This method also requires the construction of fluxes $\bq_{h,i}$ which corresponds to a larger number of degrees of freedom, but again these fluxes can be computed efficiently by solving small local problems on patches.

\begin{figure}
\includegraphics[width=0.50\textwidth]{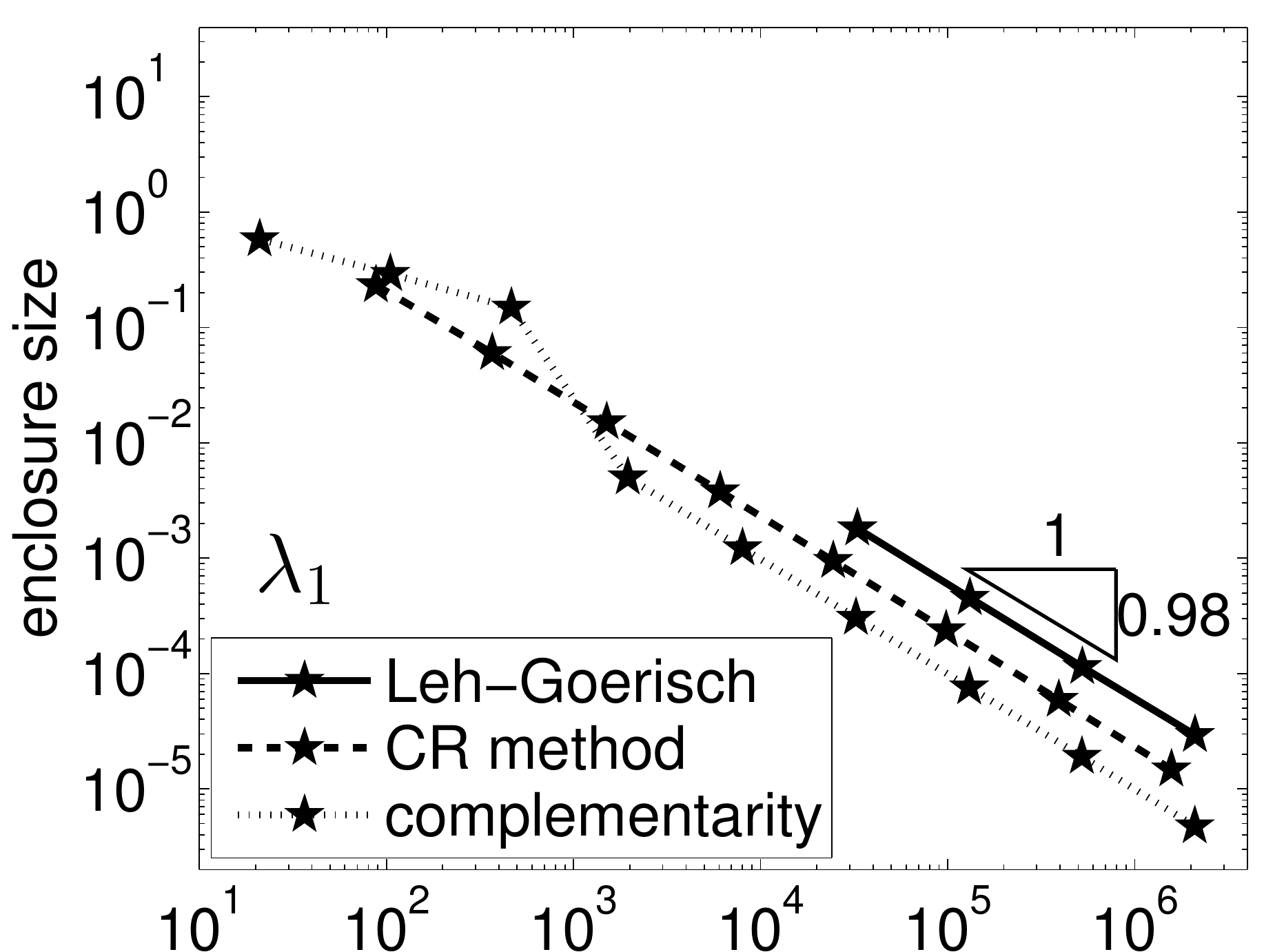}%
\includegraphics[width=0.48\textwidth]{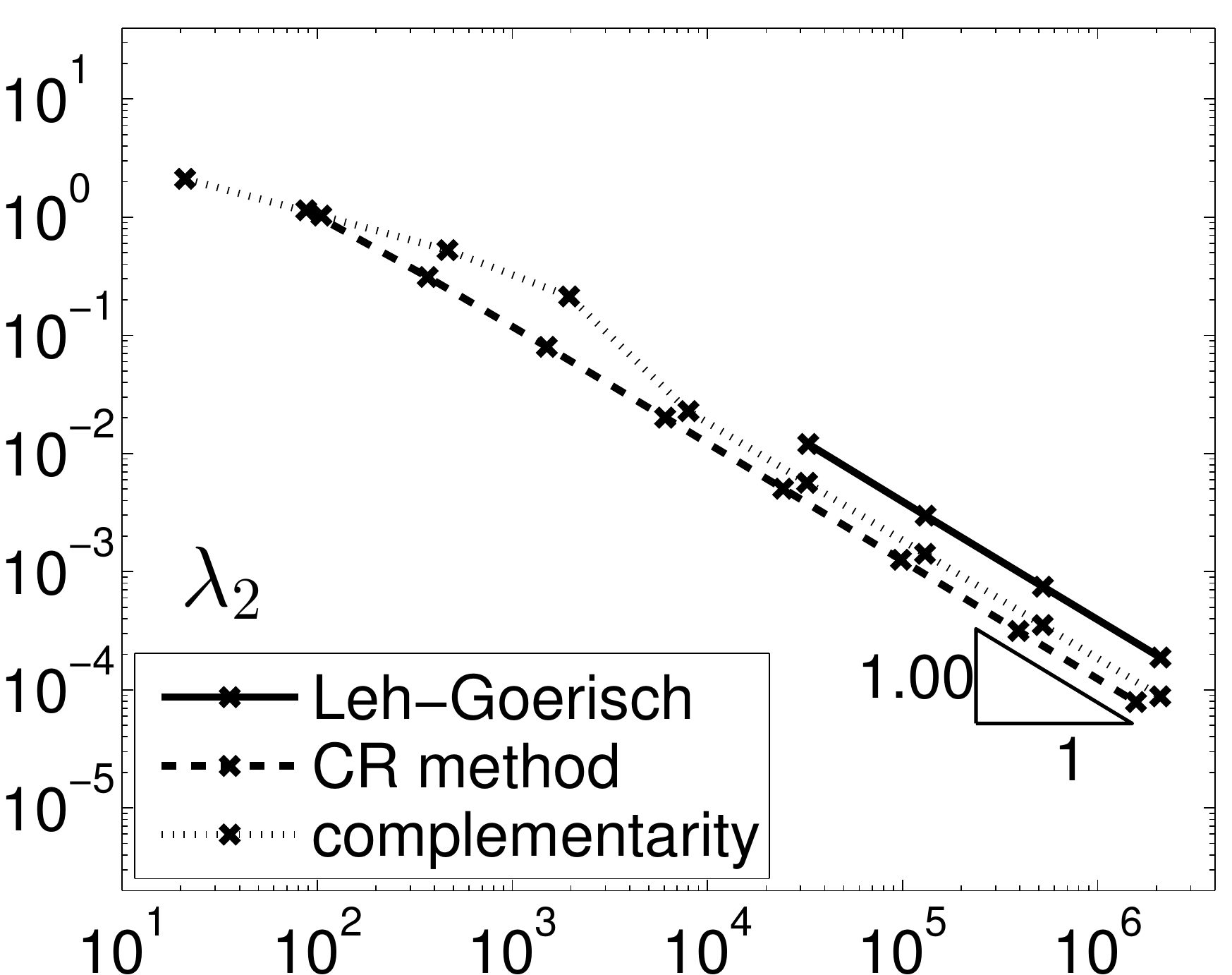}%
\\
\includegraphics[width=0.50\textwidth]{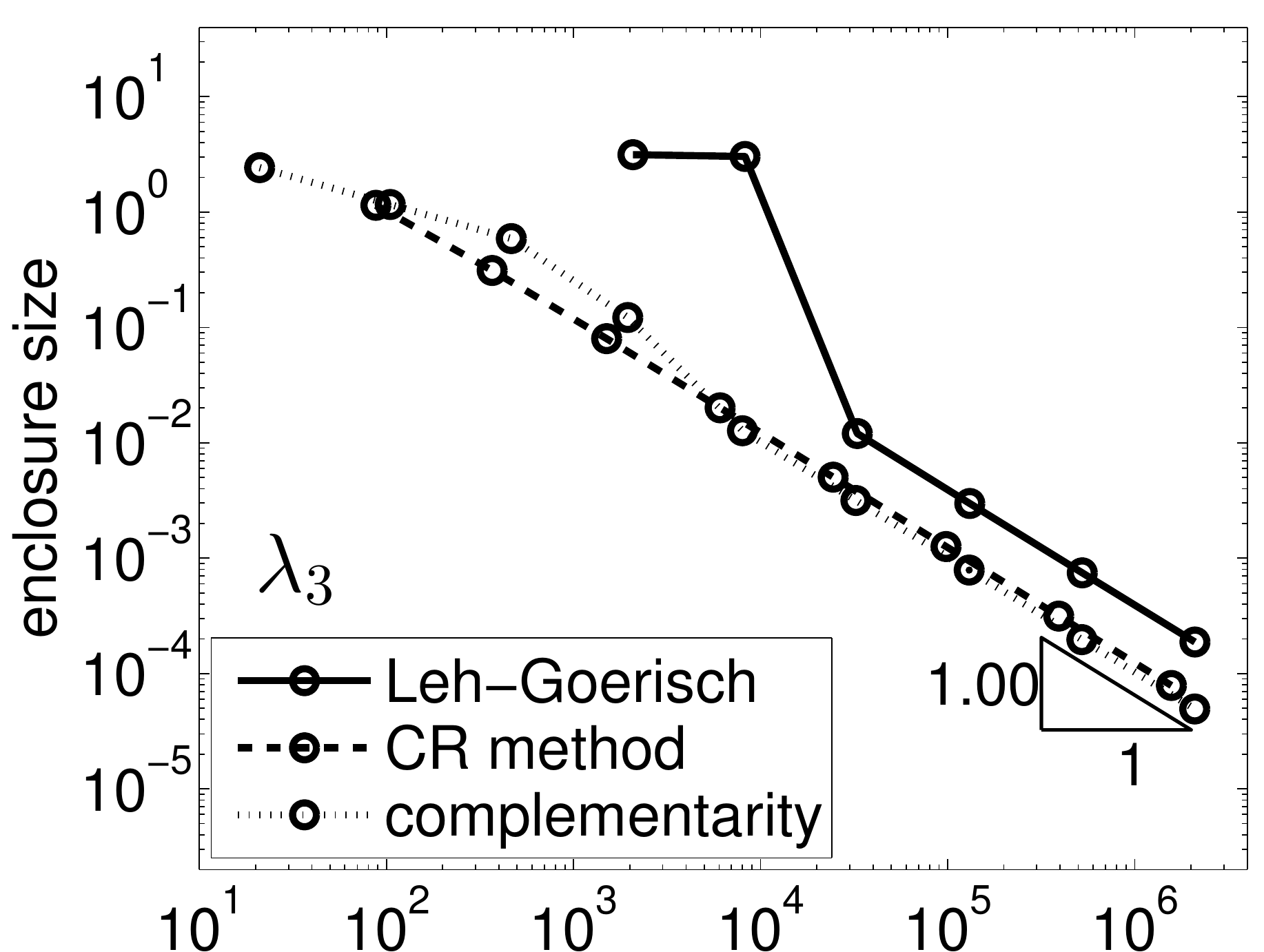}%
\includegraphics[width=0.48\textwidth]{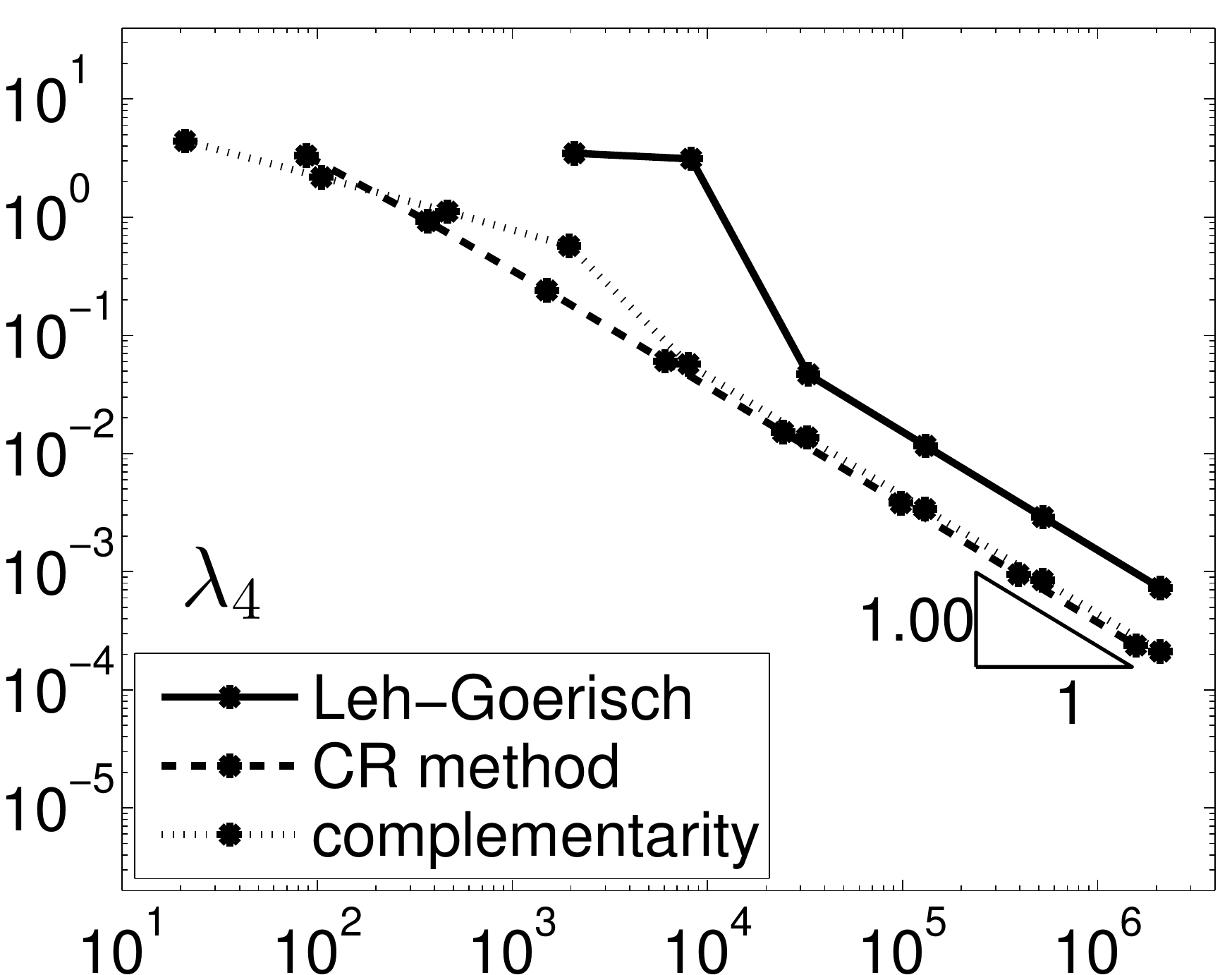}%
\\
\includegraphics[width=0.50\textwidth]{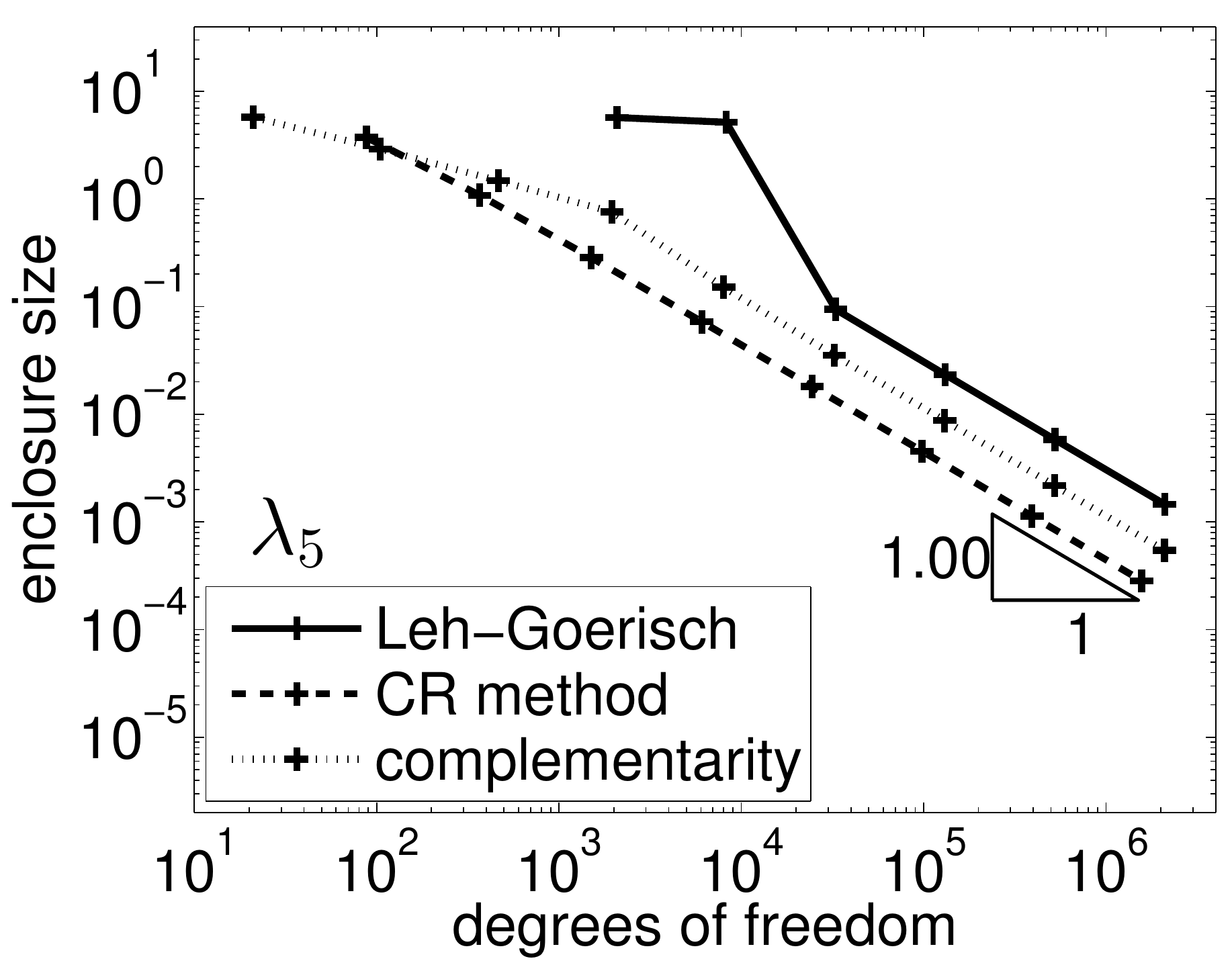}%
\includegraphics[width=0.48\textwidth]{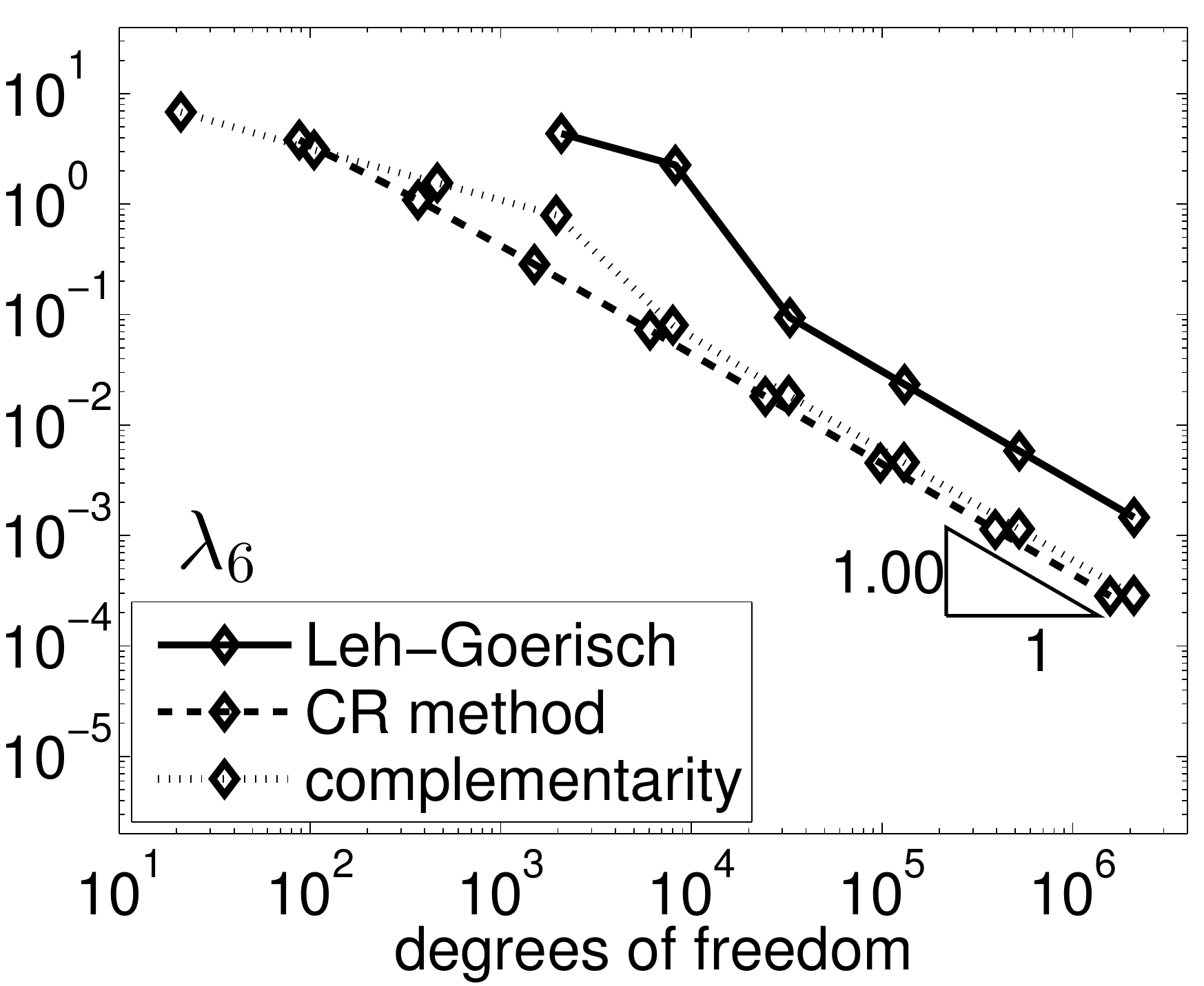}%
\caption{\label{fi:sq_encl1}
Square domain. 
Convergence of enclosure sizes of eigenvalues $\lambda_1$, $\lambda_2$, \dots, $\lambda_6$ computed by the lowest order Lehmann--Goerisch, CR, and complementarity methods.
Triangles indicate experimental orders of convergence.
}
\end{figure}

\begin{figure}
\includegraphics[width=0.50\textwidth]{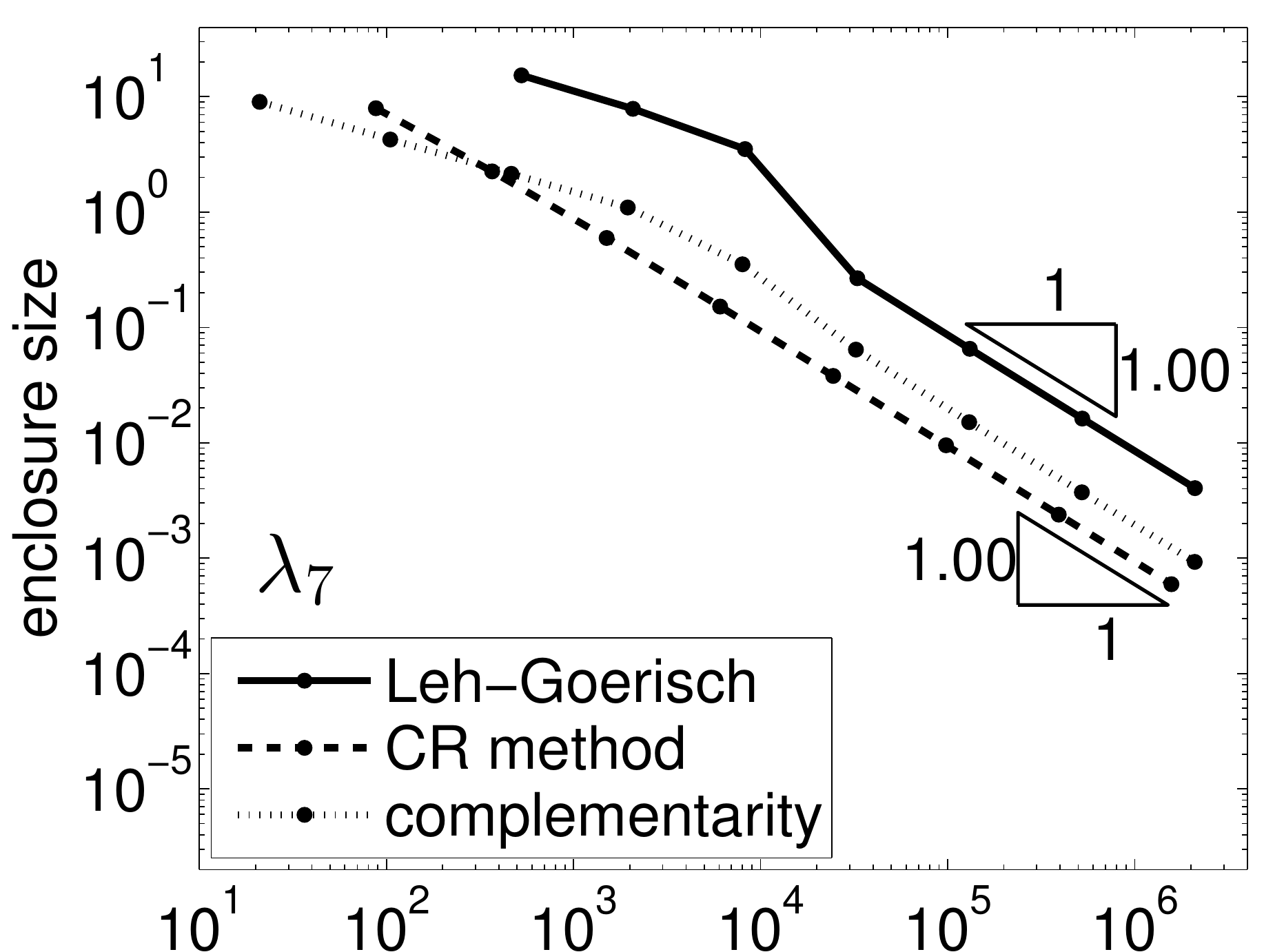}%
\includegraphics[width=0.48\textwidth]{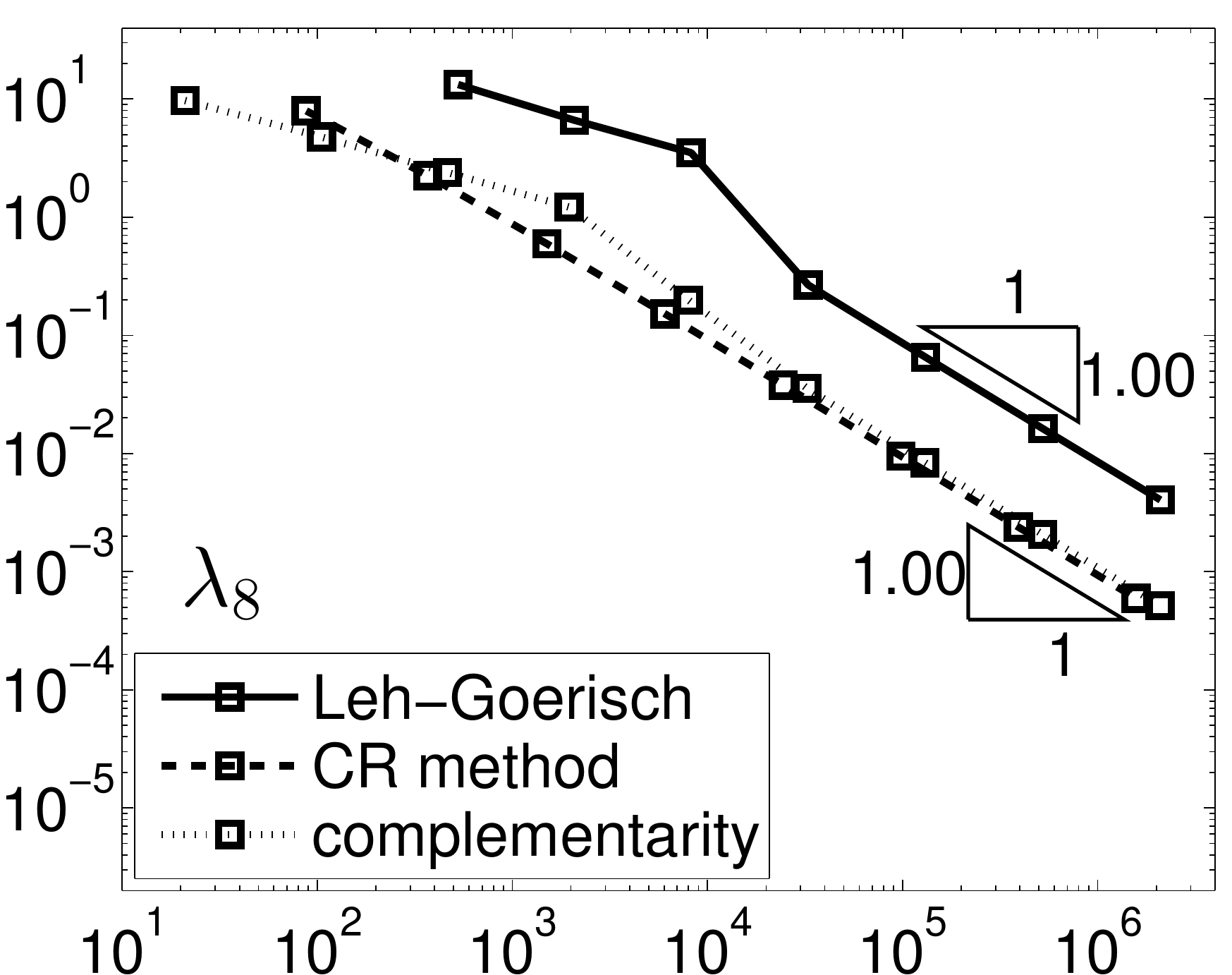}%
\\
\includegraphics[width=0.50\textwidth]{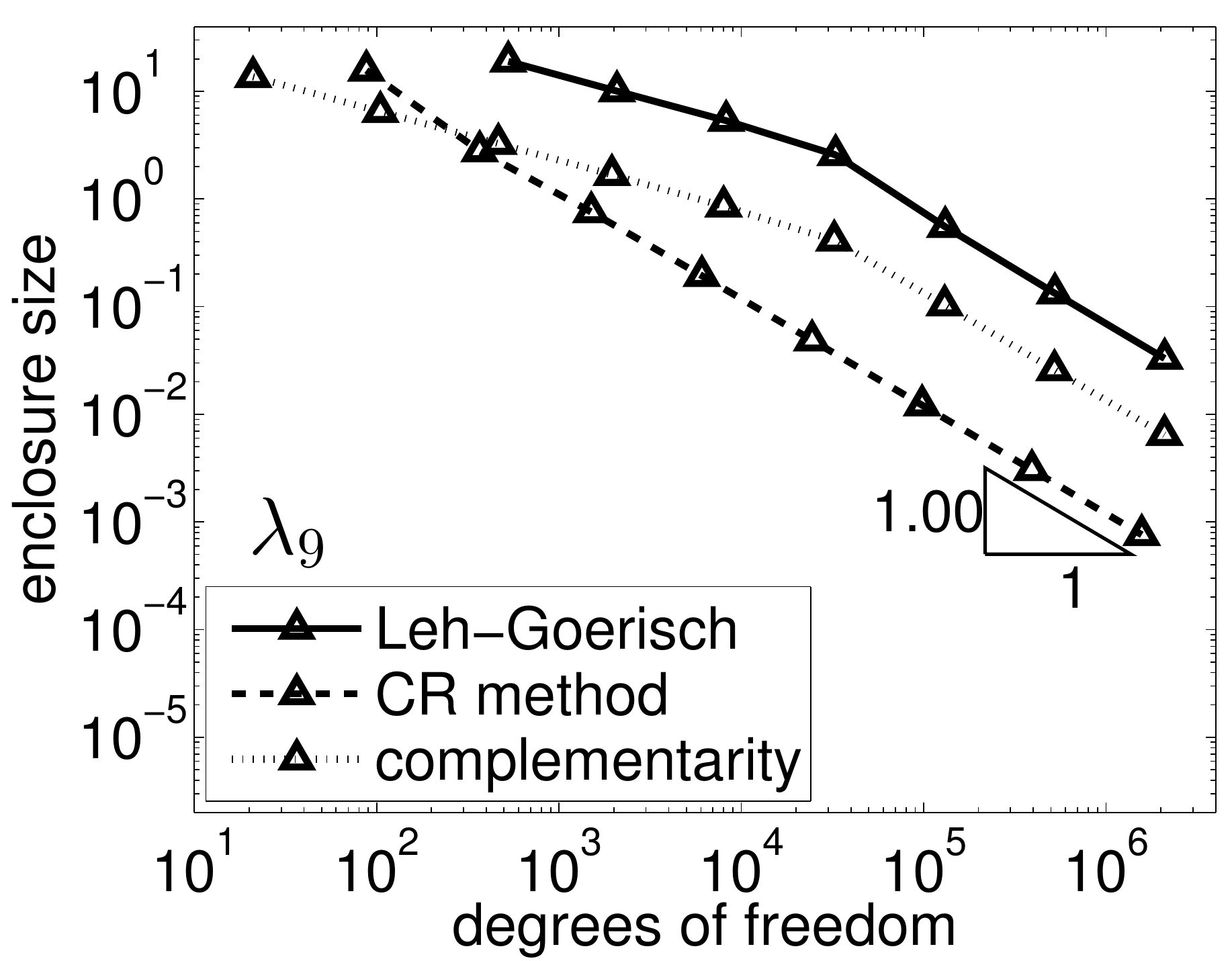}%
\includegraphics[width=0.48\textwidth]{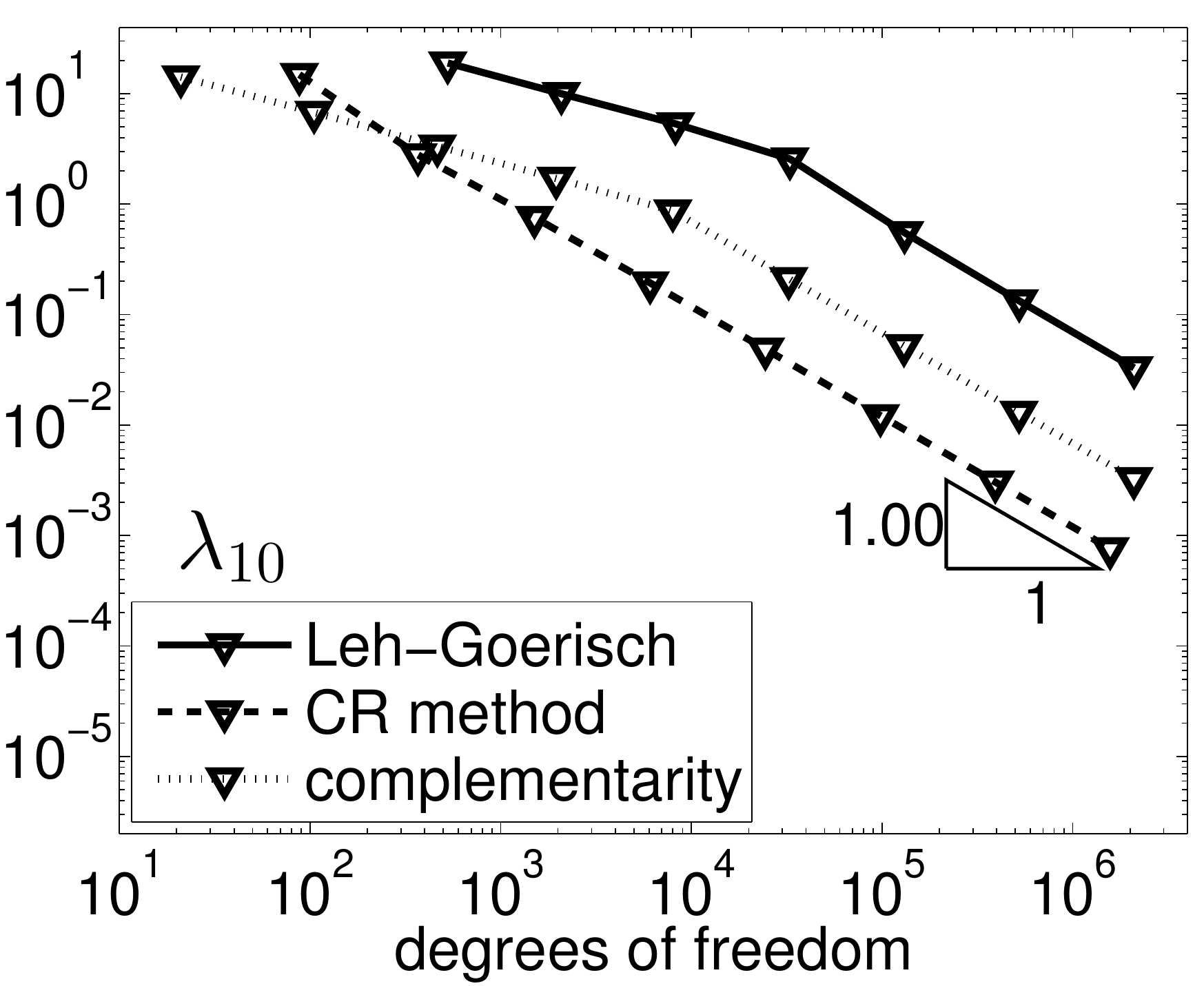}%
\caption{\label{fi:sq_encl2}
Square domain. 
Enclosure sizes of eigenvalues $\lambda_7$, $\lambda_8$, \dots, $\lambda_{10}$
obtained by the lowest order Lehmann--Goerisch, CR, and complementarity methods. 
Triangles indicate experimental orders of convergence.
}
\end{figure}

The number of degrees of freedom needed by the Lehmann--Goerisch method is the largest and we are able to refine the initial mesh 8 times. The CR method requires less degrees of freedom and allows 9 uniform refinements of the initial mesh. The complementarity method has the smallest number of degrees of freedom and enables 10 refinements.
Results obtained on the finest available meshes by the lowest order methods are reported in Table~\ref{ta:sq_eig} together with the upper bound $\Lambda_{h,i}$ computed by the Ritz--Galerkin method \eqref{eq:fem}.

Figures~\ref{fi:sq_encl1}--\ref{fi:sq_encl2} present sizes of eigenvalue enclosures, i.e. the differences between the computed upper and lower bound. More precisely the enclosure size for the eigenvalue $\lambda_i$ is given 
by $\Lambda_{h,i} - \ell_i^\incl$ for the Lehmann--Goerisch method,
by $\Lambda^*_i - \ell_i^\CR $ for the CR method, and
by $\Lambda_{h,i} - \ell_i^\cmpl$ for the complementarity method.
Figures~\ref{fi:sq_encl1}--\ref{fi:sq_encl2} show the convergence of these enclosures as
the mesh is uniformly refined and the number of degrees of freedom increases.
We immediately observe that if $N$ stands for the respective number of degrees of freedom then the enclosures converge as $O(N^{-1})$ for all three methods, see the experimental orders of convergence indicated in these figures. If $h= \max_{K\in\cT_h} \diam K$ stands for the standard mesh size parameter then this corresponds to the expected $O(h^2)$ convergence \cite{BabOsb:1991,Boffi:2010}. 

We also observe that the complementarity and CR methods provide the smallest eigenvalue enclosure with respect to the needed number of degrees of freedom. The Lehmann--Goerisch method closely follows. However, for larger eigenvalues the Lehmann--Goerisch and complementarity methods lag behind.
The Lehmann--Goerisch method requires a large number of degrees of freedom and we are not able to refine the meshes as much as for the CR and especially for the complementarity method. Therefore, we observe larger eigenvalue inclusions than for the other two methods.
On rough meshes the Lehmann--Goerisch method does not provide lower bounds on smaller eigenvalues, because the used \emph{a priori} lower bound is not sufficiently accurate. By employing an \emph{a priori} information about eigenvalues smaller than $\lambda_{11}$ we could however obtain competitive lower bounds on smaller eigenvalues even by the Lehmann--Goerisch method.

The convergence curve for the complementarity method typically shows suboptimal slope on rough meshes and optimal slope on fine meshes. This is caused by the combination of the Weinstein bound $\ell^\Wein_i$ and the Kato bound $\ell^\Kato_i$, see \eqref{eq:compl}.
We actually see $\ell^\Wein_i$ on rough meshes and $\ell^\Kato_i$ on fine meshes.


Further observation is that enclosure sizes of lower eigenvalues are smaller than
enclosure sizes of higher eigenvalues. 
All methods provide less accurate results as the index of the eigenvalue increases. 
This effect is highlighted in our results, because we present absolute sizes of enclosures. However, we would observe the loss of accuracy for the higher eigenvalues even if we plotted relative sizes of enclosures weighted by sizes of corresponding eigenvalues.

\section{Numerical results -- a dumbbell shaped domain}
\label{se:dumbbell}

This section applies the three methods described in Sections~\ref{se:2}--\ref{se:4} 
to compute lower and upper bounds on the first $m=10$ eigenvalues of the problem proposed in \cite{TreBet2006}. It is the eigenvalue problem \eqref{eq:modprob} posed in the dumbbell shaped domain $\Omega = (0,\pi)^2 \cup [\pi,5\pi/4]\times(3\pi/8,5\pi/8) \cup (5\pi/4,9\pi/4)\times(0,\pi)$, see Figure~\ref{fi:sq_initmesh} (middle). 
This is a more realistic example, where the exact eigenvalues are not known.
For consistency with the CR method, the Lehmann--Goerisch and complementarity methods are considered of order $k=1$.

\begin{table}
\centerline{
\begin{tabular}{ccccc}
 & lower bound       & lower bound & lower bound     & upper bound    \\ 
 & Lehmann--Goerisch & CR method   & complementarity & Ritz--Galerkin  \\ \hline
$\lambda_{1}$  &   1.9556896  &   1.95565230 &   1.95569083 &   1.95582583 \\ 
$\lambda_{2}$  &   1.96058965 &   1.96055131 &   1.96064783 &   1.96071159 \\ 
$\lambda_{3}$  &   4.79930938 &   4.80005018 &   4.79979095 &   4.80091560 \\ 
$\lambda_{4}$  &   4.82860260 &   4.82924108 &   4.82940402 &   4.83002932 \\ 
$\lambda_{5}$  &   4.99635717 &   4.99648459 &   4.99667320 &   4.99686964 \\ 
$\lambda_{6}$  &   4.99637117 &   4.99649844 &   4.99678524 &   4.99688342 \\ 
$\lambda_{7}$  &   7.98319666 &   7.98599709 &   7.98582092 &   7.98704483 \\ 
$\lambda_{8}$  &   7.98326415 &   7.98606433 &   7.98650019 &   7.98711174 \\ 
$\lambda_{9}$  &   9.34945616 &   9.35480997 &   9.33633170 &   9.35772093 \\ 
$\lambda_{10}$ &   9.50362132 &   9.50864166 &   9.49999525 &   9.51119420 \\ 
\end{tabular}
}
\caption{\label{ta:dmbl_eig}
Dumbbell shaped domain. Eigenvalue bounds computed on the finest meshes. The mesh sizes of the finest meshes were $h=0.0368$ for the Lehmann--Goerisch method, $h=0.0184$ for the CR method, and $h=0.0092$ for the complementarity method.
}
\end{table}

\begin{figure}
\includegraphics[width=0.50\textwidth]{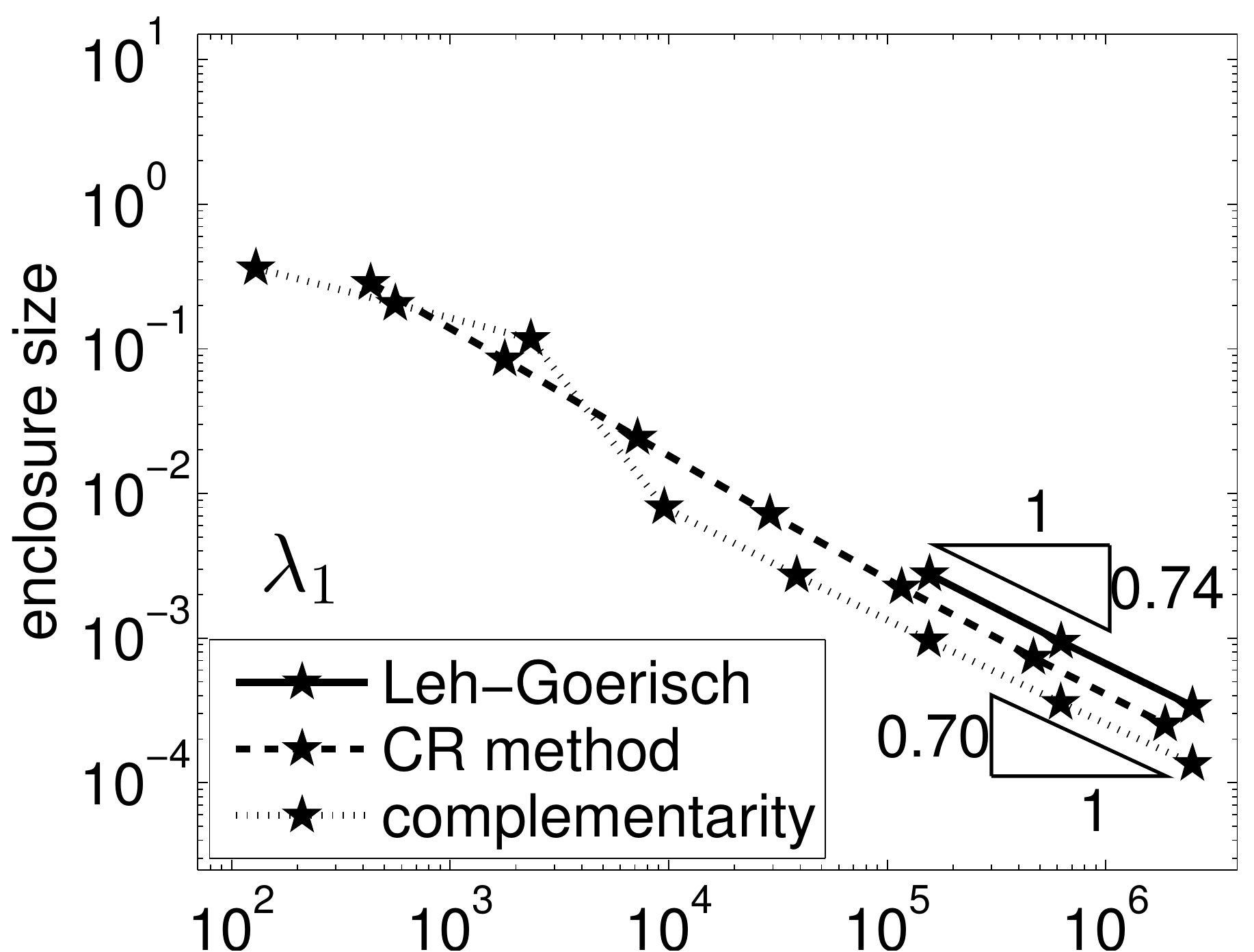}%
\includegraphics[width=0.48\textwidth]{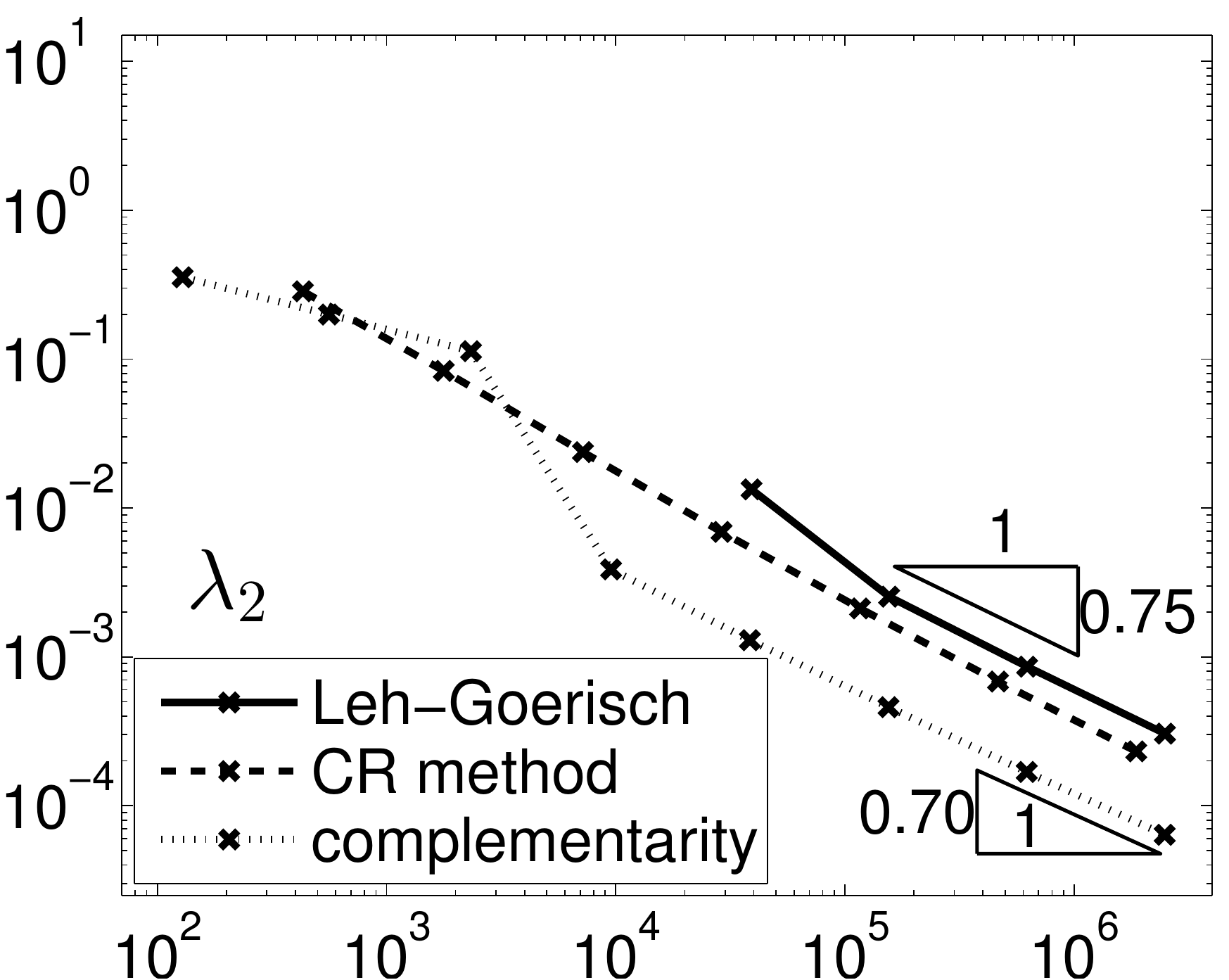}%
\\
\includegraphics[width=0.50\textwidth]{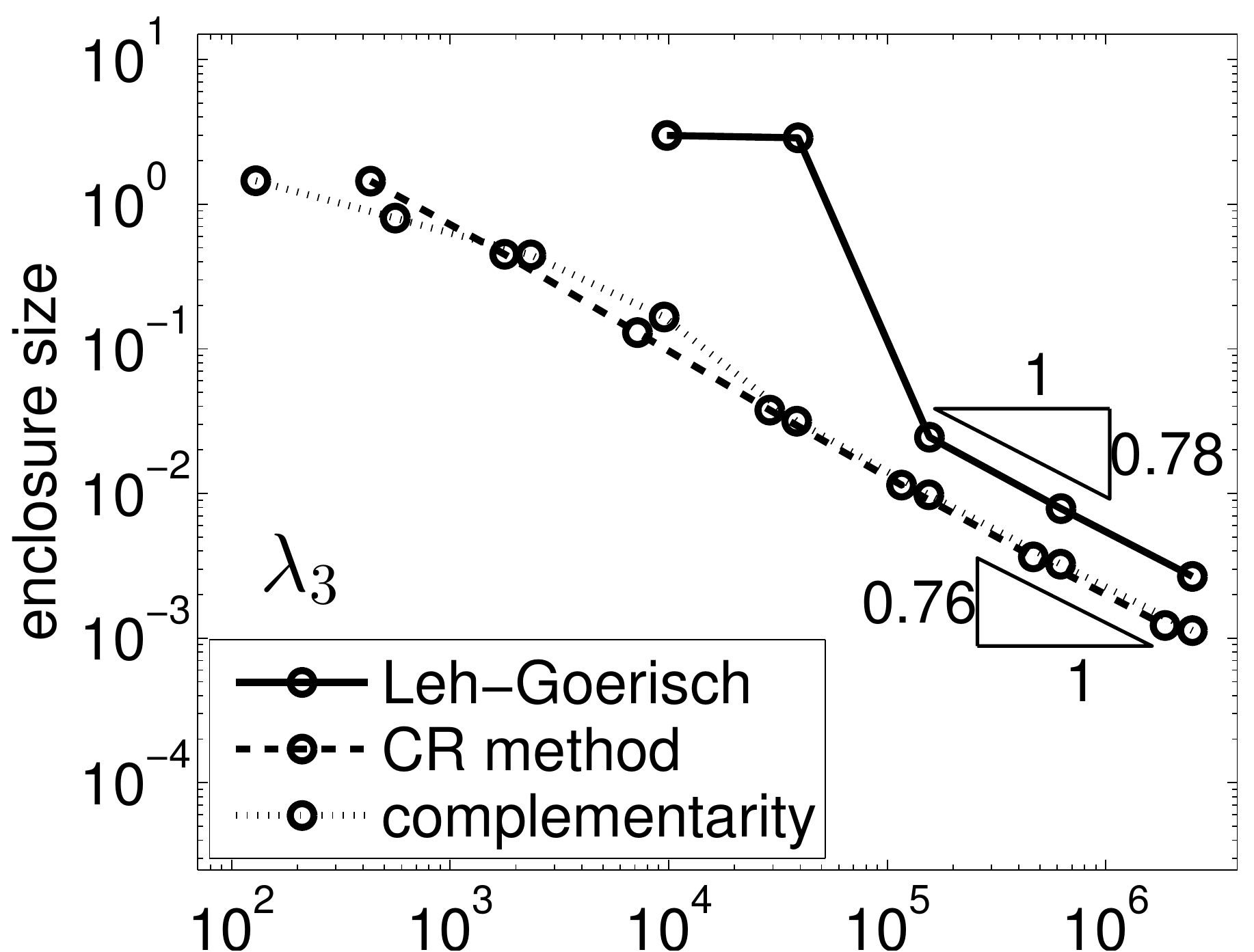}%
\includegraphics[width=0.48\textwidth]{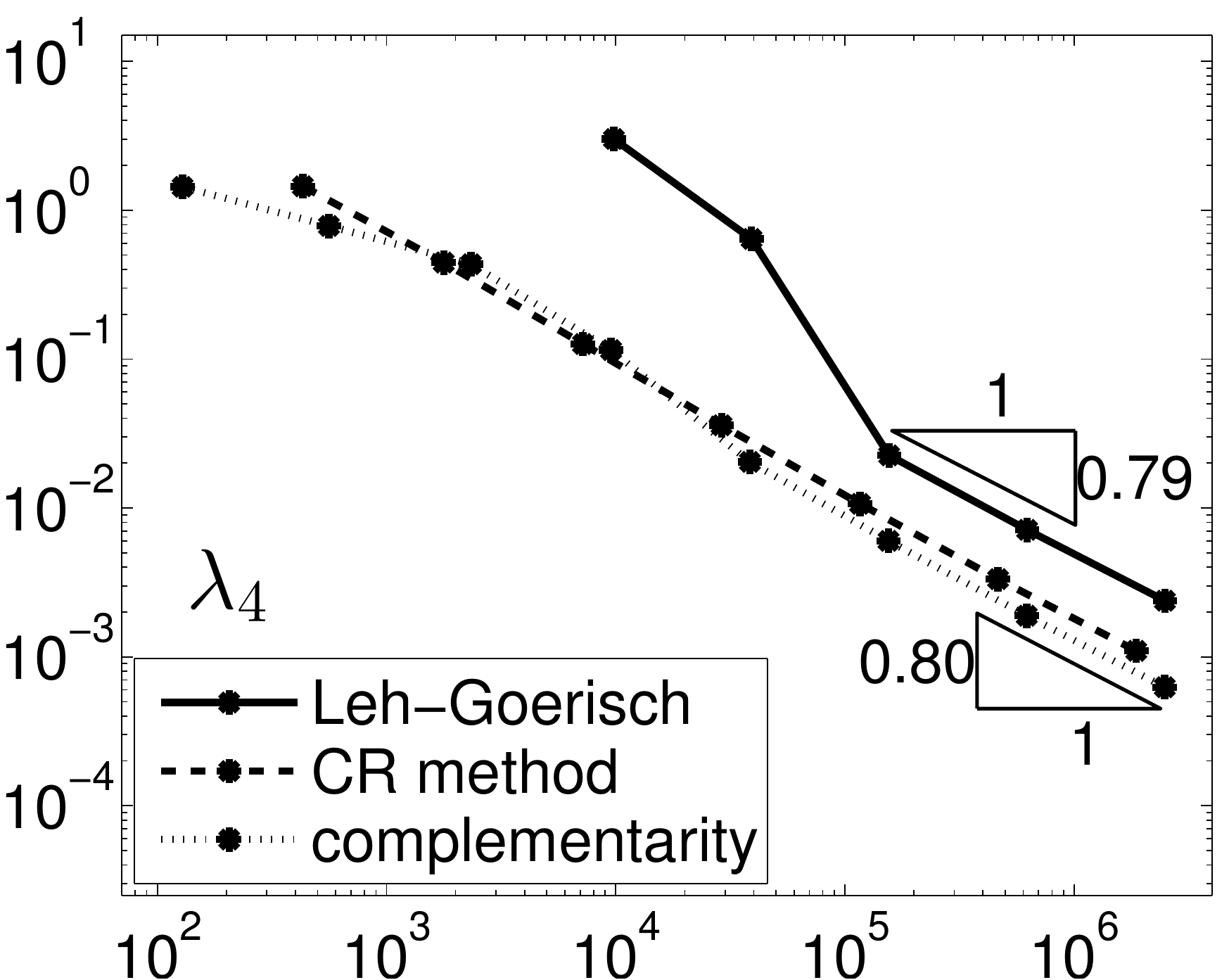}%
\\
\includegraphics[width=0.50\textwidth]{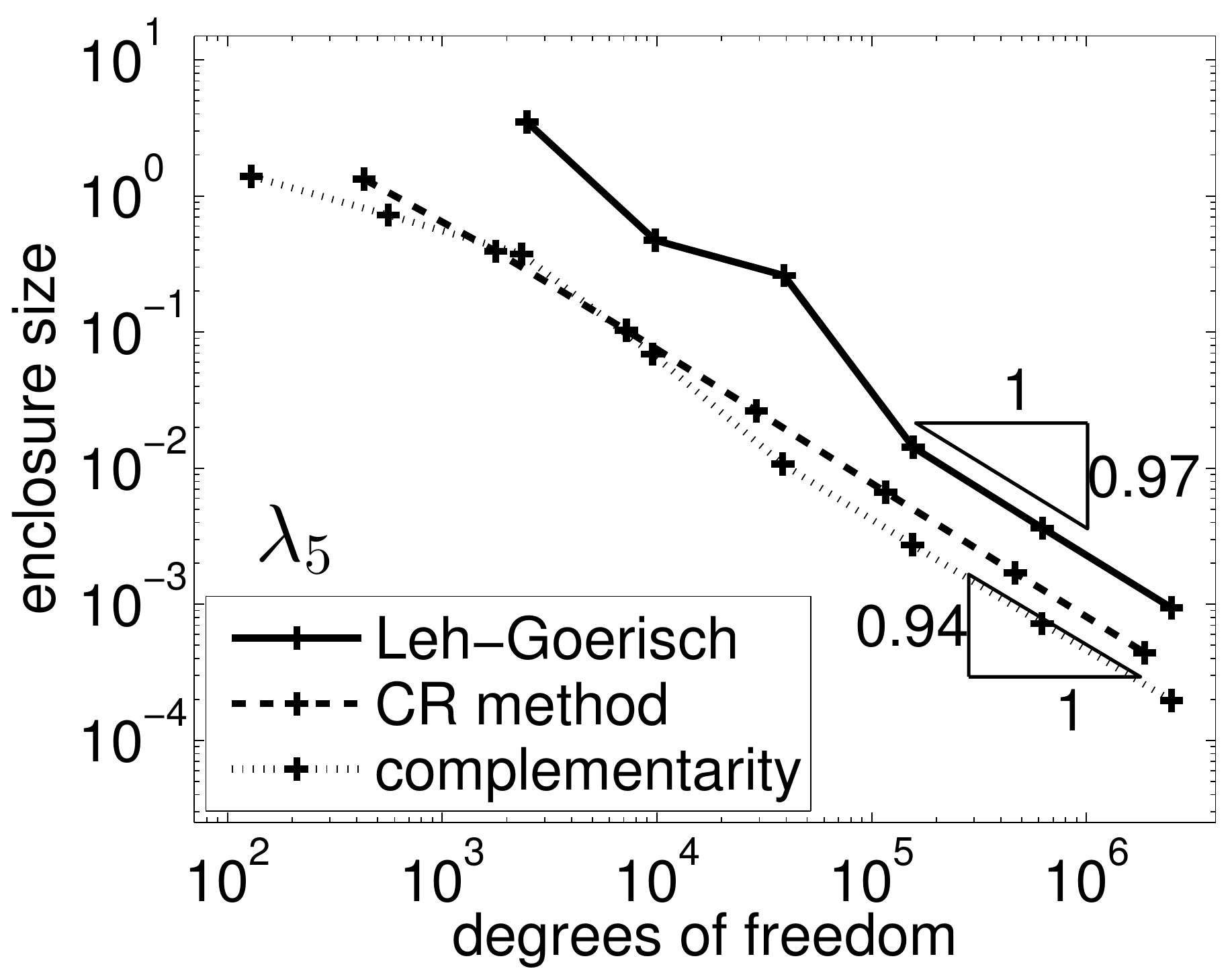}%
\includegraphics[width=0.48\textwidth]{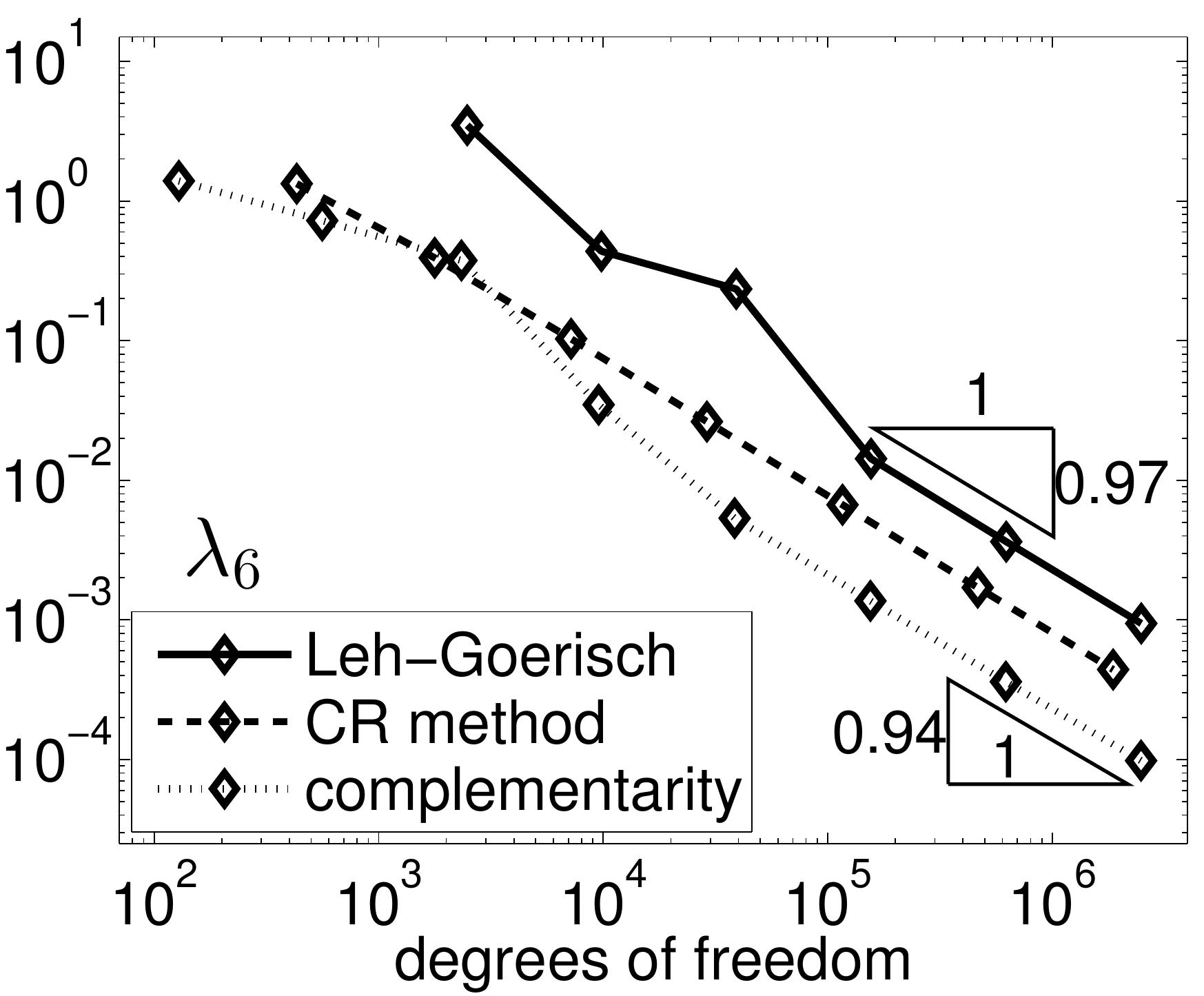}%
\caption{\label{fi:dmbl_encl1}
Dumbbell shaped domain.
Enclosure sizes of eigenvalues $\lambda_1$, $\lambda_2$, \dots, $\lambda_6$
obtained by  the lowest order Lehmann--Goerisch, CR, and complementarity methods. 
Triangles indicate experimental orders of convergence. 
}
\end{figure}

\begin{figure}
\includegraphics[width=0.50\textwidth]{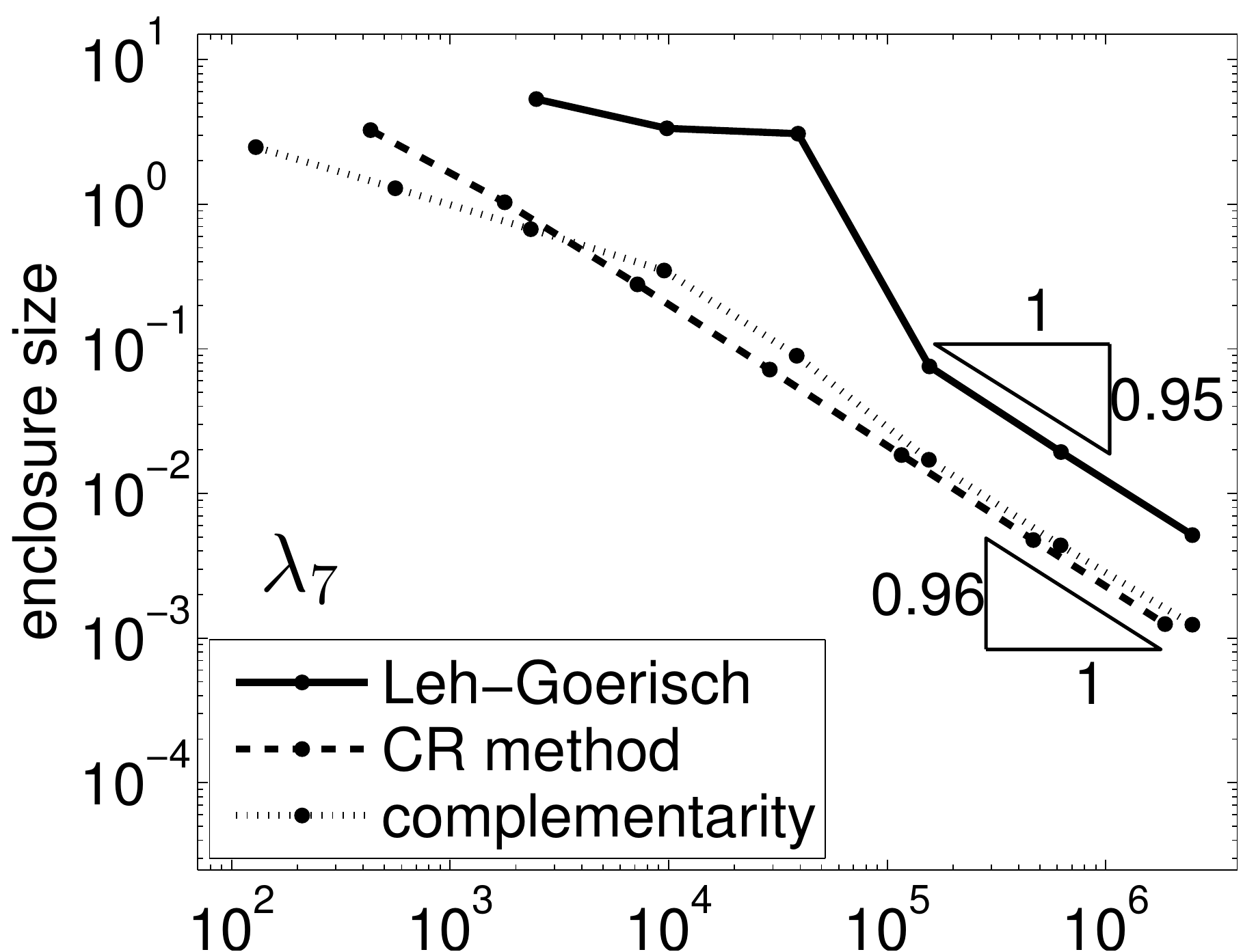}%
\includegraphics[width=0.48\textwidth]{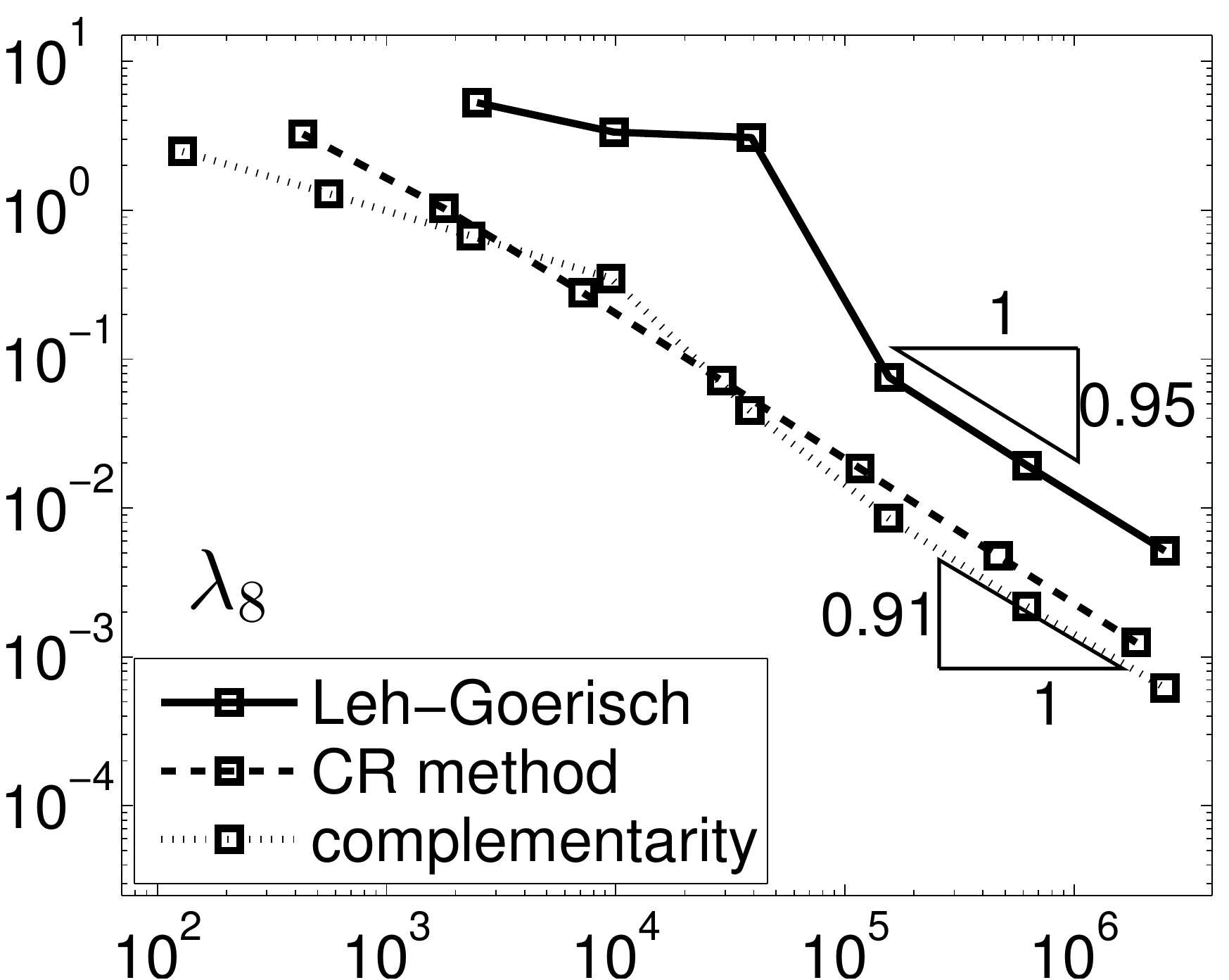}%
\\
\includegraphics[width=0.50\textwidth]{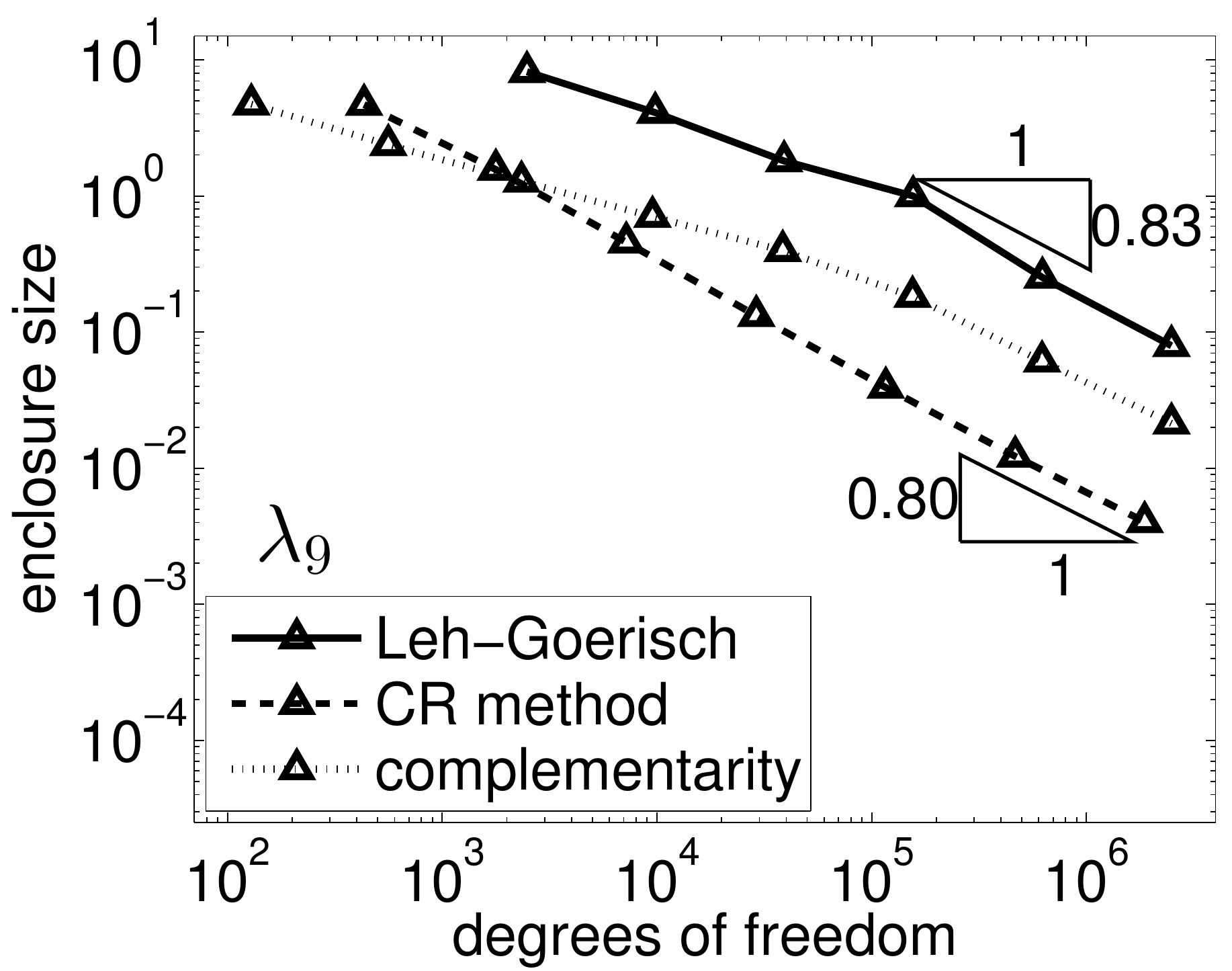}%
\includegraphics[width=0.48\textwidth]{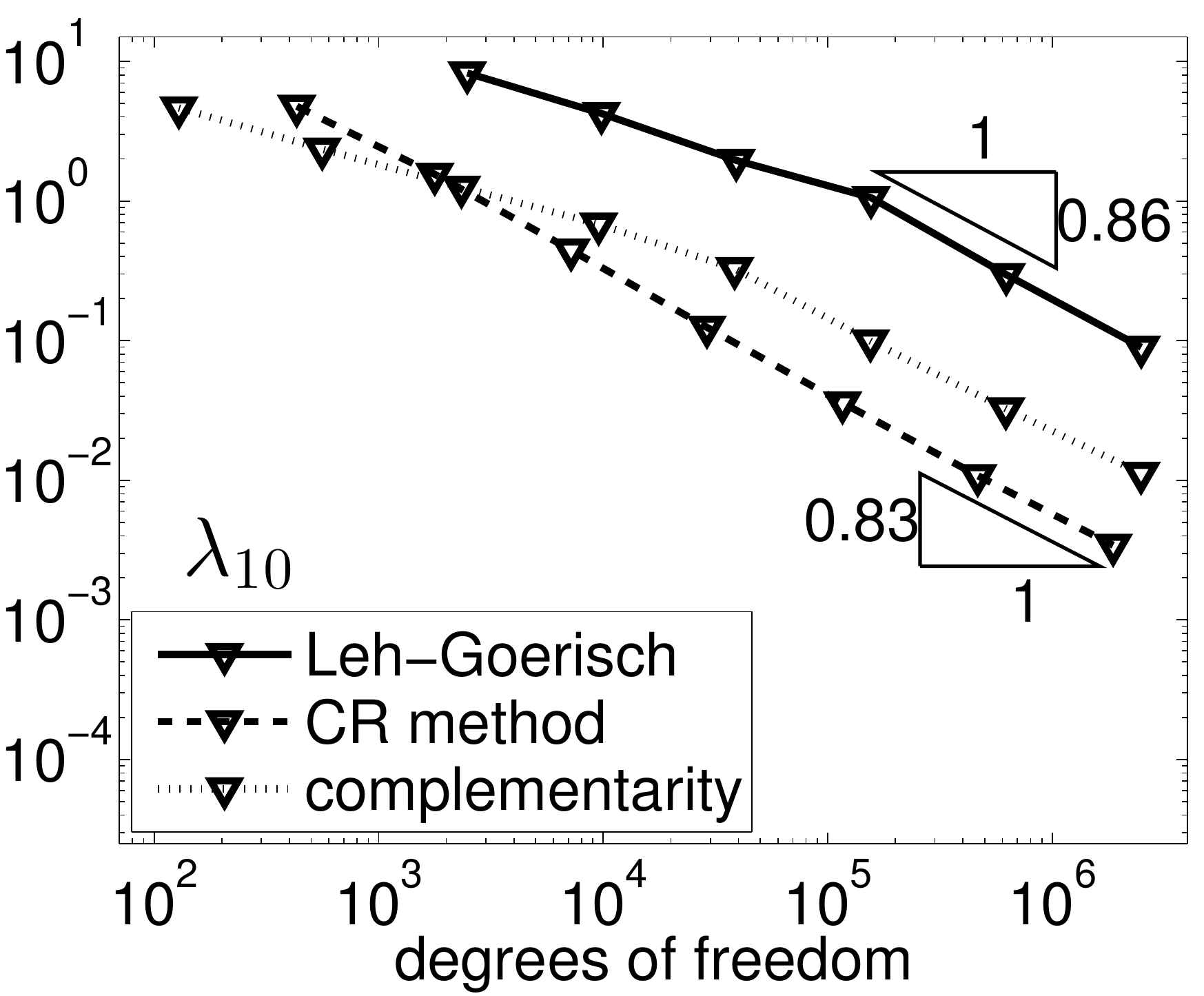}%
\caption{\label{fi:dmbl_encl2}
Dumbbell shaped domain.
Enclosure sizes of eigenvalues $\lambda_7$, $\lambda_8$, \dots, $\lambda_{10}$
obtained by  the lowest order Lehmann--Goerisch, CR, and complementarity methods. 
Triangles indicate experimental orders of convergence. 
}
\end{figure}

The three methods are compared using the same methodology as in the previous section. 
The initial triangulation is depicted in Figure~\ref{fi:sq_initmesh} (middle).
For the Lehmann--Goerisch method we were able to refine the mesh 6 times, for the CR method 7 times, and for the complementarity method 8 times.
Table~\ref{ta:dmbl_eig} reports the computed lower and upper bounds on the finest available meshes
for all three methods.
We observe that although the two-sided bounds are quite tight, the enclosing intervals overlap for $\lambda_5$ and $\lambda_6$ and for $\lambda_7$ and $\lambda_8$. On this level of accuracy, we cannot decide whether these two pairs are multiple or isolated eigenvalues. On the other hand, $\lambda_1$ and $\lambda_2$ as well as $\lambda_3$ and $\lambda_4$ are tight pairs of eigenvalues and the computed two-sided bounds are sufficiently accurate to show that they are all isolated.

Figures~\ref{fi:dmbl_encl1}--\ref{fi:dmbl_encl2} present sizes of eigenvalue enclosures for eigenvalues $\lambda_1$, $\lambda_2$, \dots, $\lambda_{10}$ and their dependence on the number of degrees freedom. We recall that the number of degrees of freedom corresponds to
$N^\mathrm{mix} = \dim \bWh + \dim Q_h$ for the Lehmann--Goerisch method, 
to $N^\CR = \dim \VCR_h$ for the CR method, and to $N^\mathrm{conf} = \dim V_h$ for the complementarity method.
In this example the complementarity and CR methods provide again the most accurate results with respect to the needed number of degrees of freedom. The convergence rates are in several cases spoiled by singularities of eigenvectors at reentrant corners of the domain. Optimal rates could be achieved by the adaptive mesh refinement, but this is not the goal of this paper. 
As in the previous example, the Lehmann--Goerisch method yields only slightly less accurate results than the other two methods for the smallest eigenvalues, however, for larger eigenvalues the differences are larger. The wider enclosure sizes of the Lehmann--Goerisch method are due to the large number of degrees of freedom it requires, because we are not able to refine the meshes as much as for the CR and especially the complementary method.
In accordance with the previous example absolute sizes of eigenvalue enclosures increase for higher eigenvalues for all three methods.


\section{Numerical results -- chopped off square}
\label{se:ctsqr}

In order to test the presented methods on a domain without obvious symmetries, we consider the domain $\Omega$ depicted in Figure~\ref{fi:sq_initmesh} (right). It is the square $(0,\pi)^2$ with an asymmetrically removed corner. The removed corner is the right-angle triangle with vertices $(0.8\pi,\pi)$, $(\pi,0.9\pi)$, and $(\pi,\pi)$.
Again, lower and upper bounds are computed by the three presented methods for the first $m=10$ eigenvalues.
The order of the Lehmann--Goerisch and complementarity methods is $k=1$ for consistency.
The exact eigenvalues are not known, but they are supposed to be close to eigenvalues of the square computed in Section~\ref{se:sqrpipi}.

\begin{table}
\centerline{
\begin{tabular}{ccccc}
 & lower bound       & lower bound & lower bound     & upper bound    \\ 
 & Lehmann--Goerisch & CR method   & complementarity & Ritz--Galerkin \\ \hline
$\lambda_{1}$ &   2.00428875 &   2.00428589 &   2.00429114 &   2.00429304 \\ 
$\lambda_{2}$ &   4.99999175 &   4.99999709 &   5.00000218 &   5.00004190 \\ 
$\lambda_{3}$ &   5.03006005 &   5.03006713 &   5.03009172 &   5.03011214 \\ 
$\lambda_{4}$ &   8.05214209 &   8.05226769 &   8.05228797 &   8.05238463 \\ 
$\lambda_{5}$ &   9.99994012 &   10.0003563 &   10.0002639 &   10.0005314 \\ 
$\lambda_{6}$ &   10.0547637 &   10.0551768 &   10.0552239 &   10.0553595 \\ 
$\lambda_{7}$ &   12.9989536 &   13.0004808 &   13.0003594 &   13.0007870 \\ 
$\lambda_{8}$ &   13.1960062 &   13.1977998 &   13.1978755 &   13.1981073 \\ 
$\lambda_{9}$ &   16.9839761 &   17.0022290 &   16.9998586 &   17.0027346 \\ 
$\lambda_{10}$&   17.0432010 &   17.0643490 &   17.0633160 &   17.0648692 \\ 
\end{tabular}
}
\caption{\label{ta:ctsqr_eig}
Chopped off square. Eigenvalue bounds computed on the finest meshes. The mesh sizes of the finest meshes were $h=0.0112$ for the Lehmann--Goerisch method, $h=0.0056$ for the CR method, and $h=0.0028$ for the complementarity method.
}
\end{table}

\begin{figure}
\includegraphics[width=0.50\textwidth]{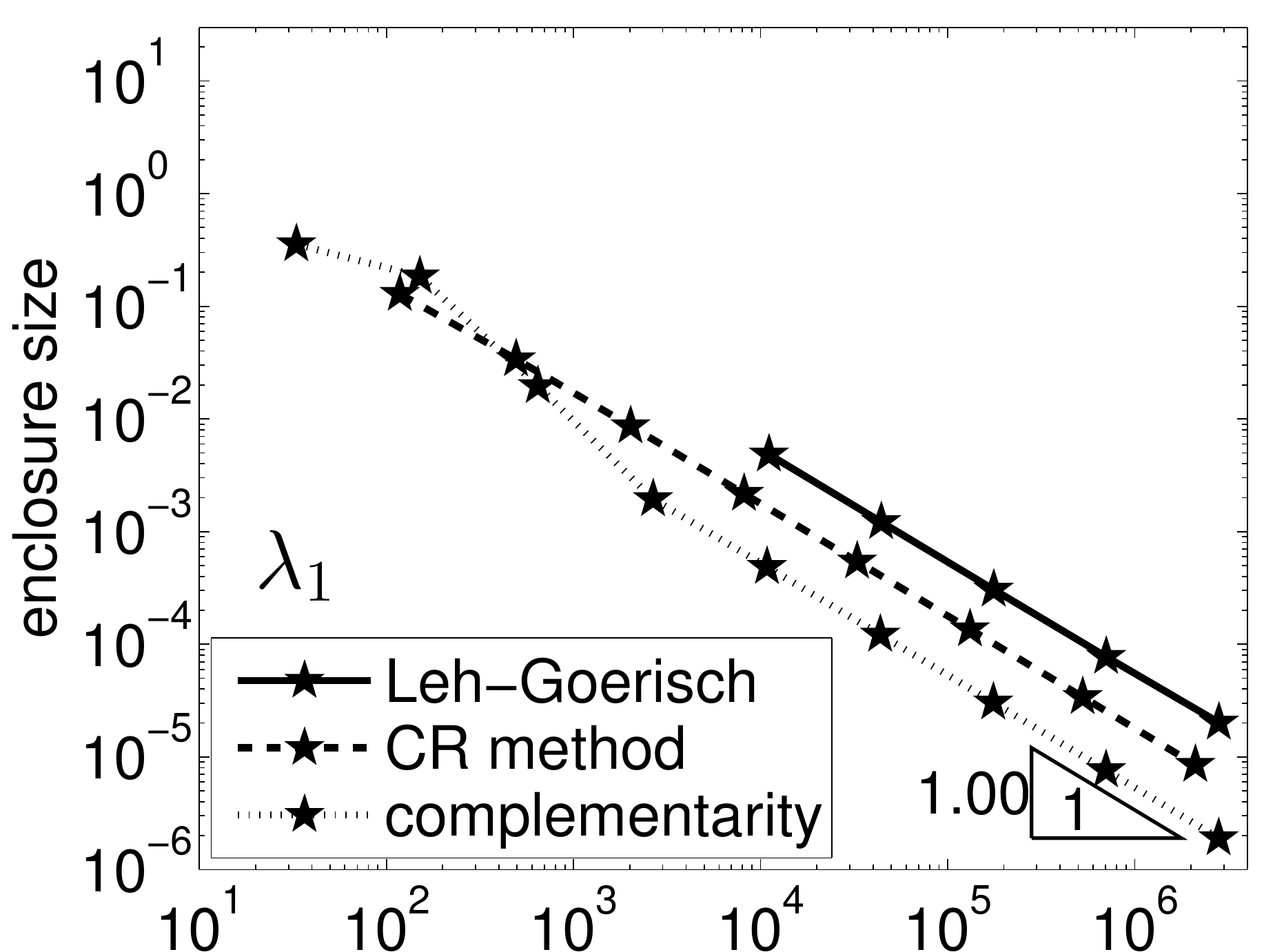}%
\includegraphics[width=0.48\textwidth]{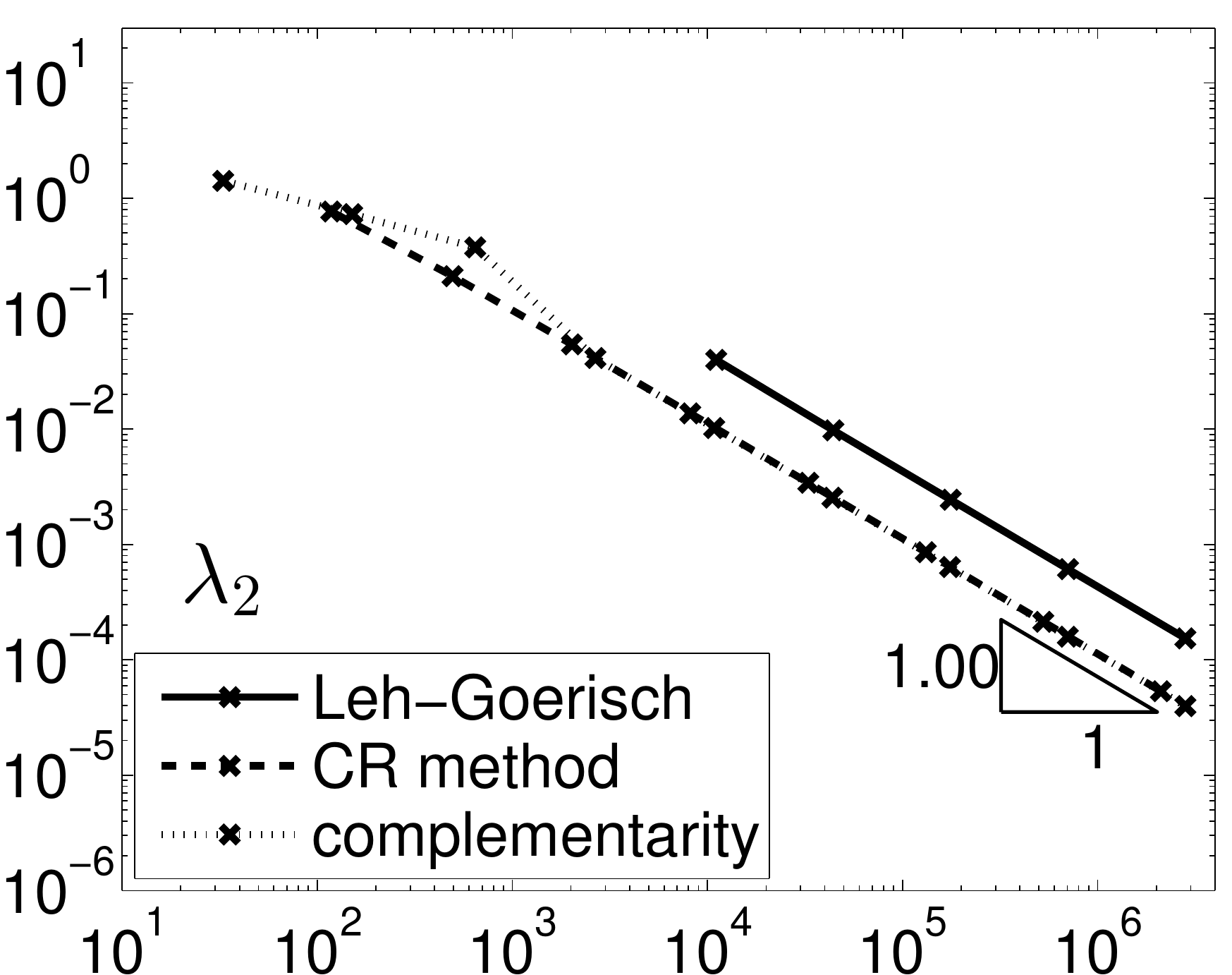}%
\\
\includegraphics[width=0.50\textwidth]{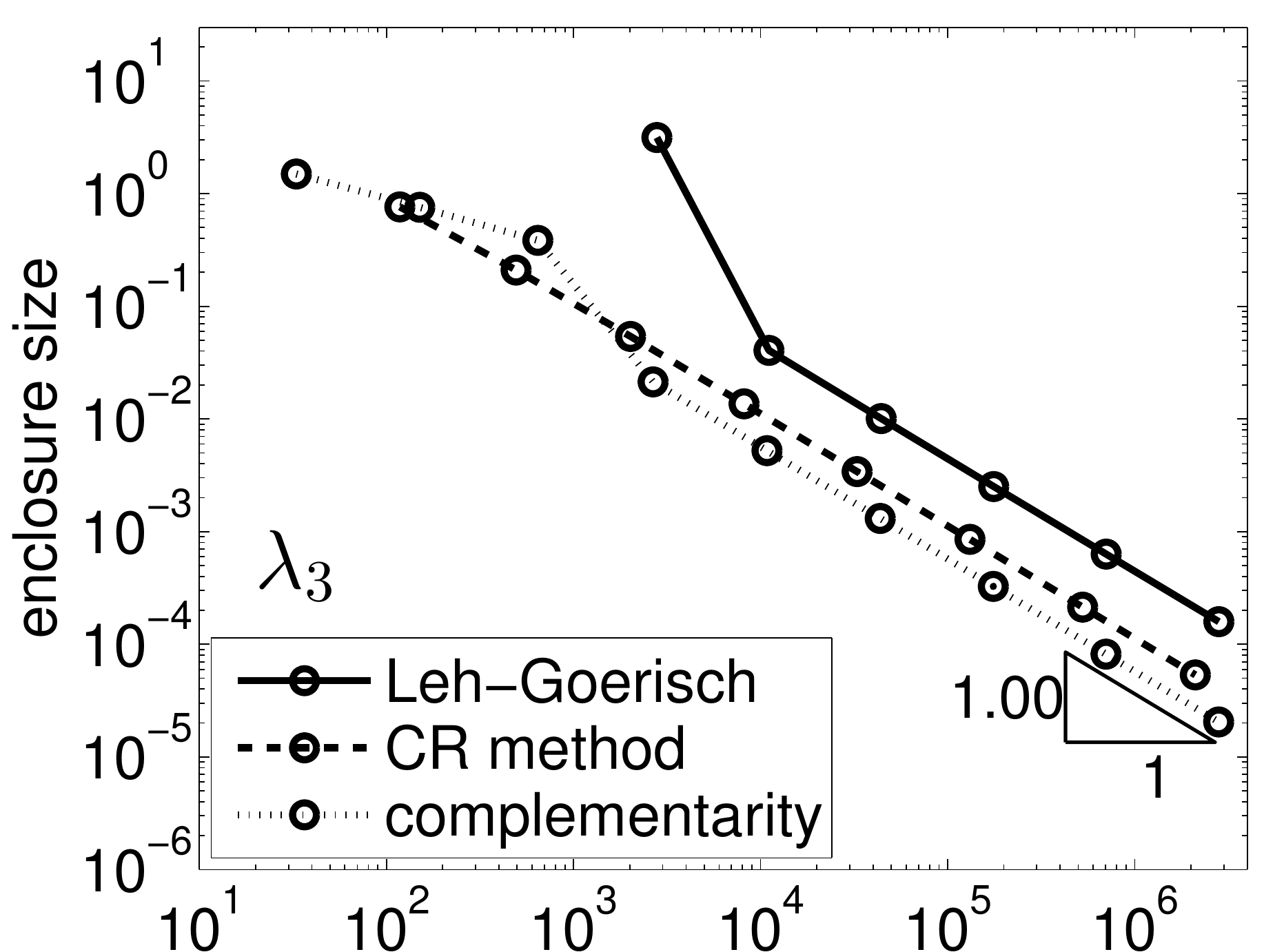}%
\includegraphics[width=0.48\textwidth]{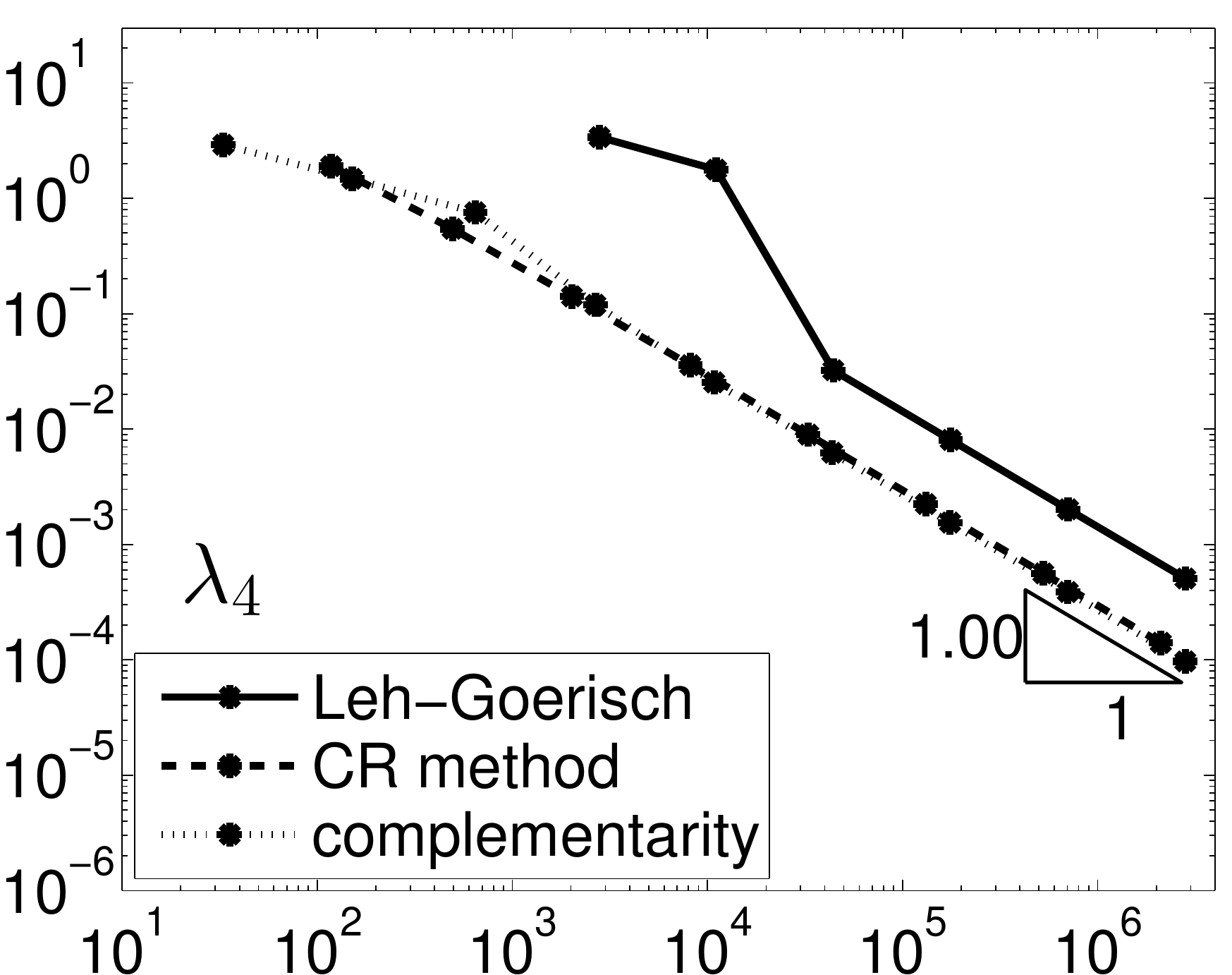}%
\\
\includegraphics[width=0.50\textwidth]{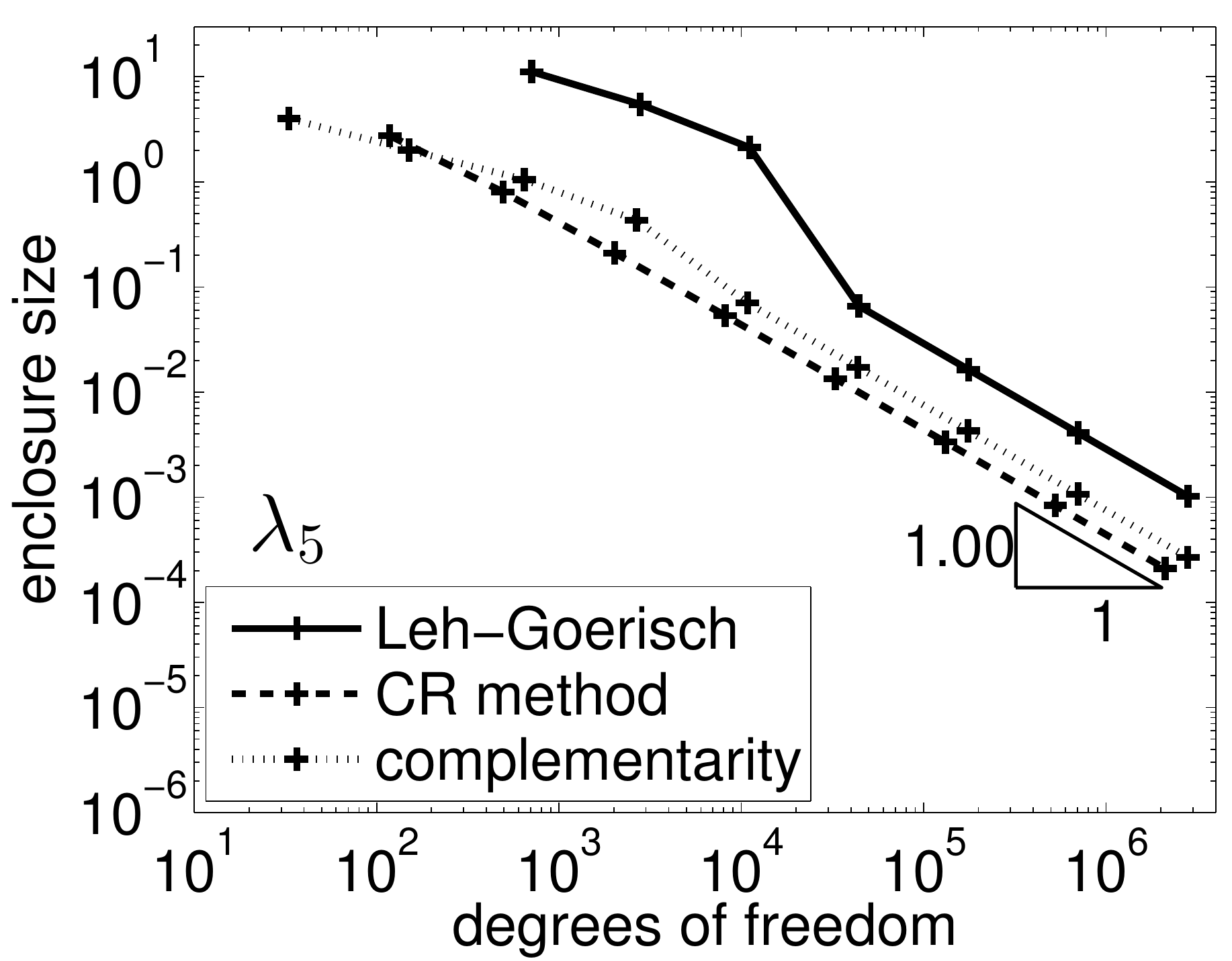}%
\includegraphics[width=0.48\textwidth]{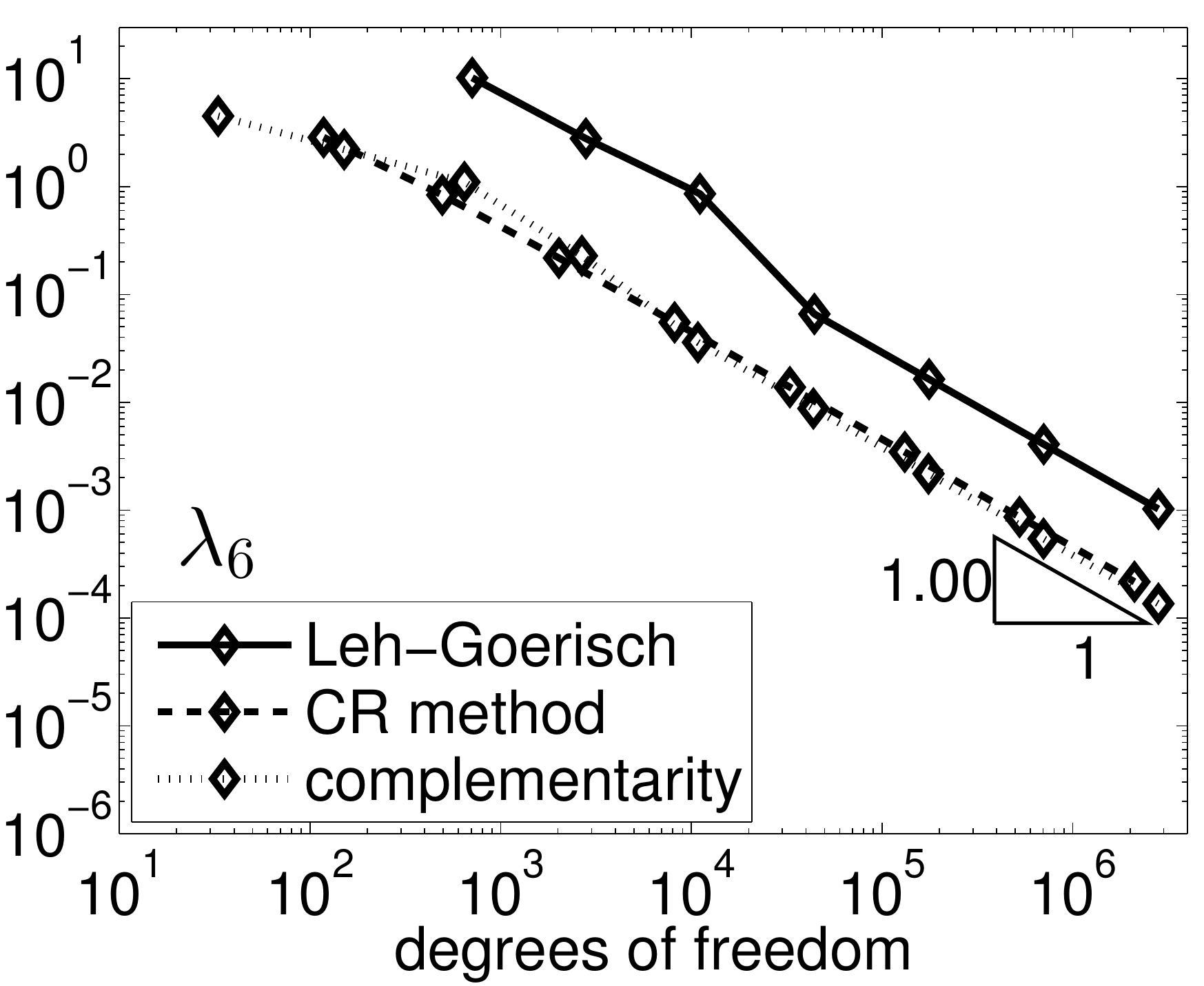}%

\caption{\label{fi:ctsqr_encl1}
Chopped off square.
Enclosure sizes of eigenvalues $\lambda_1$, $\lambda_2$, \dots, $\lambda_6$
obtained by the lowest order Lehmann--Goerisch, CR, and complementarity methods. 
Triangles indicate experimental orders of convergence. 
}
\end{figure}

\begin{figure}
\includegraphics[width=0.50\textwidth]{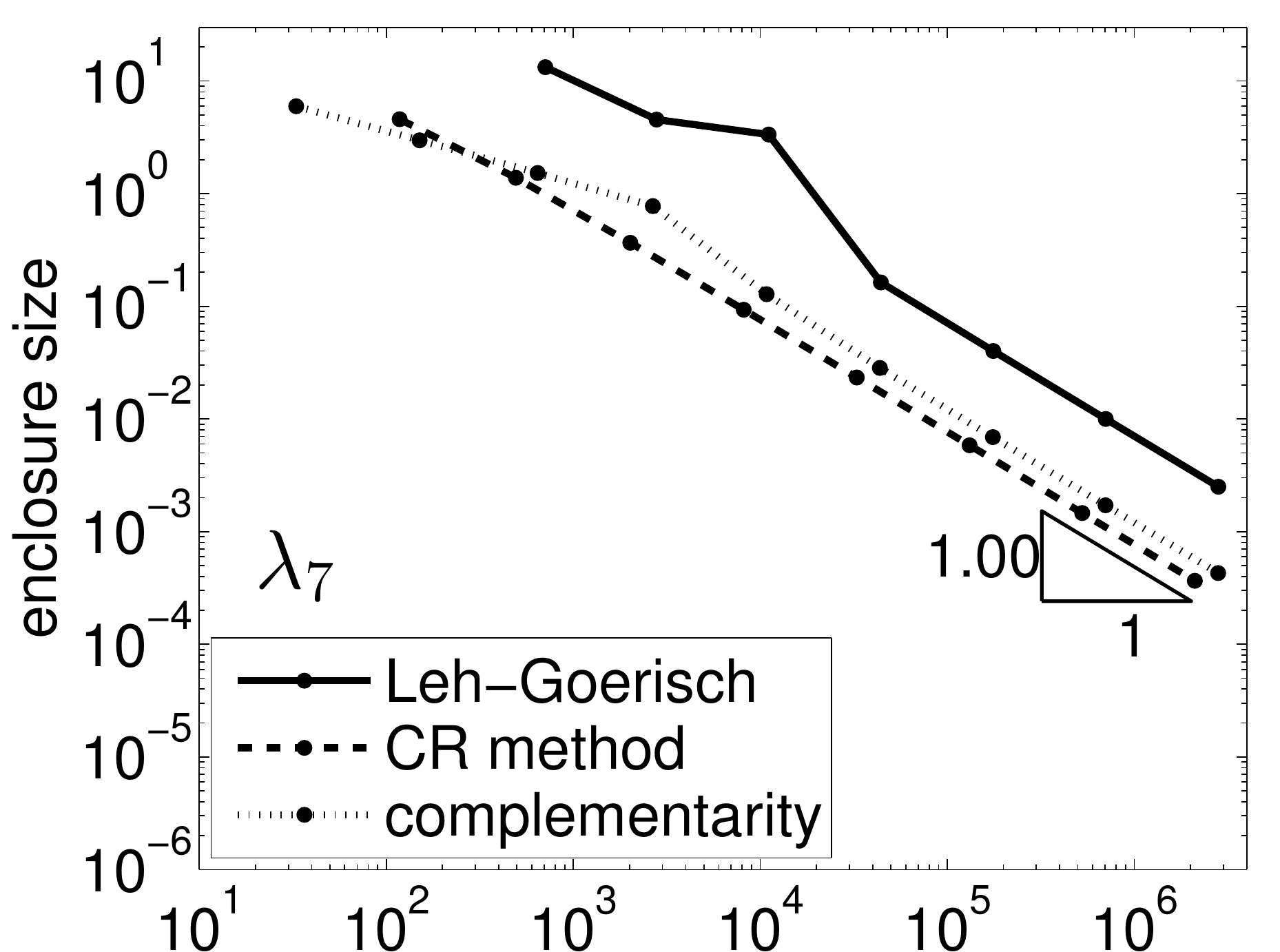}%
\includegraphics[width=0.48\textwidth]{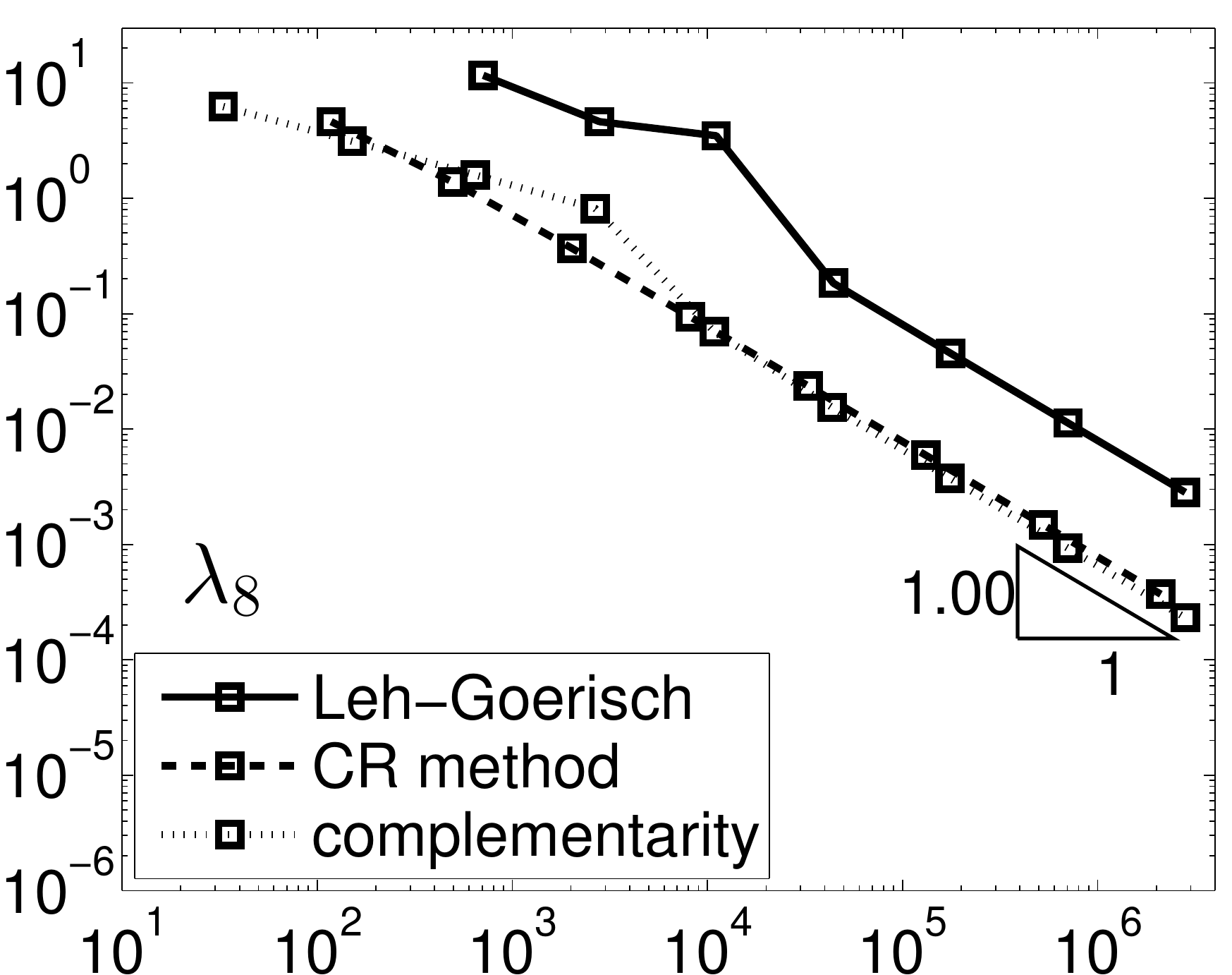}%
\\
\includegraphics[width=0.50\textwidth]{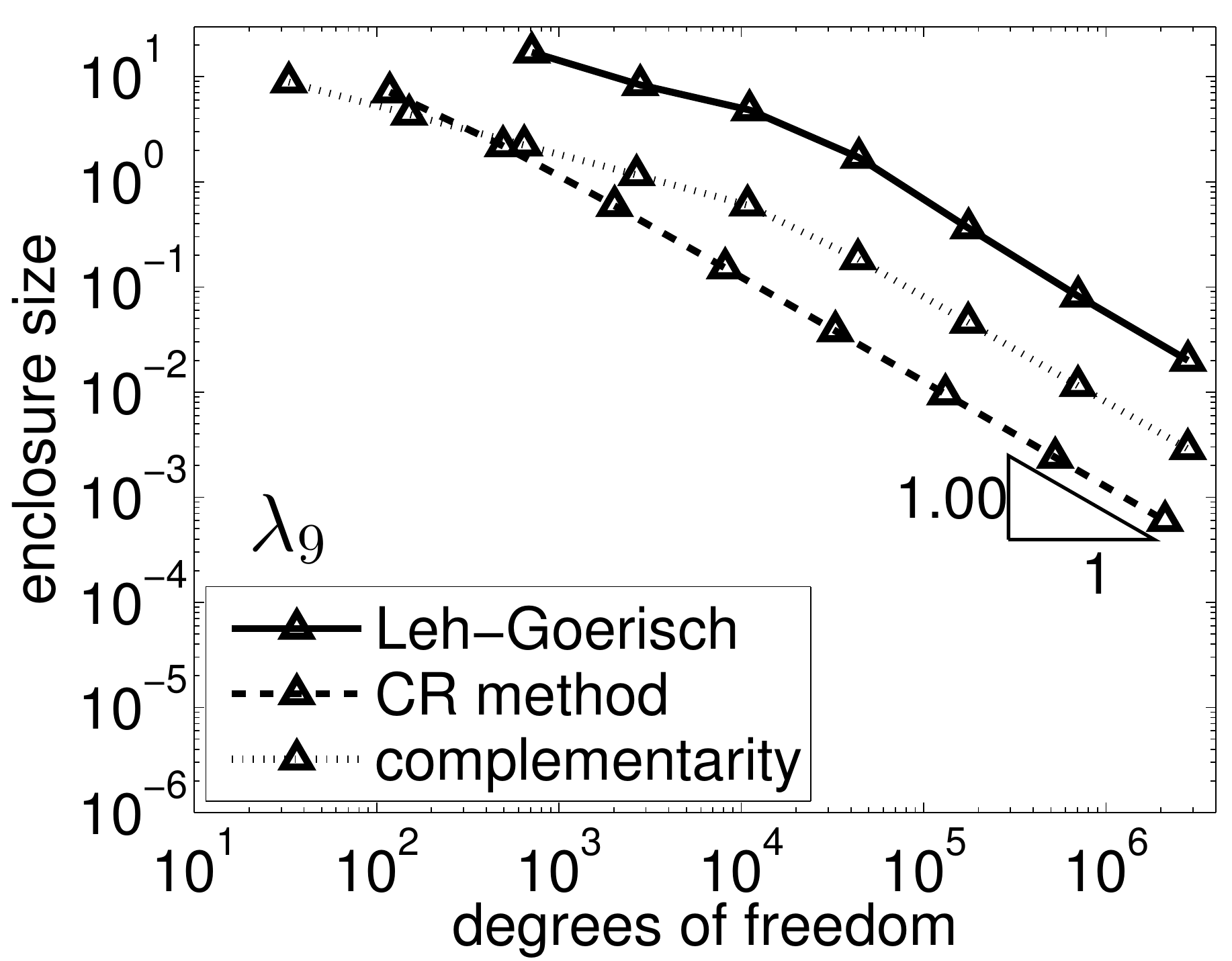}%
\includegraphics[width=0.48\textwidth]{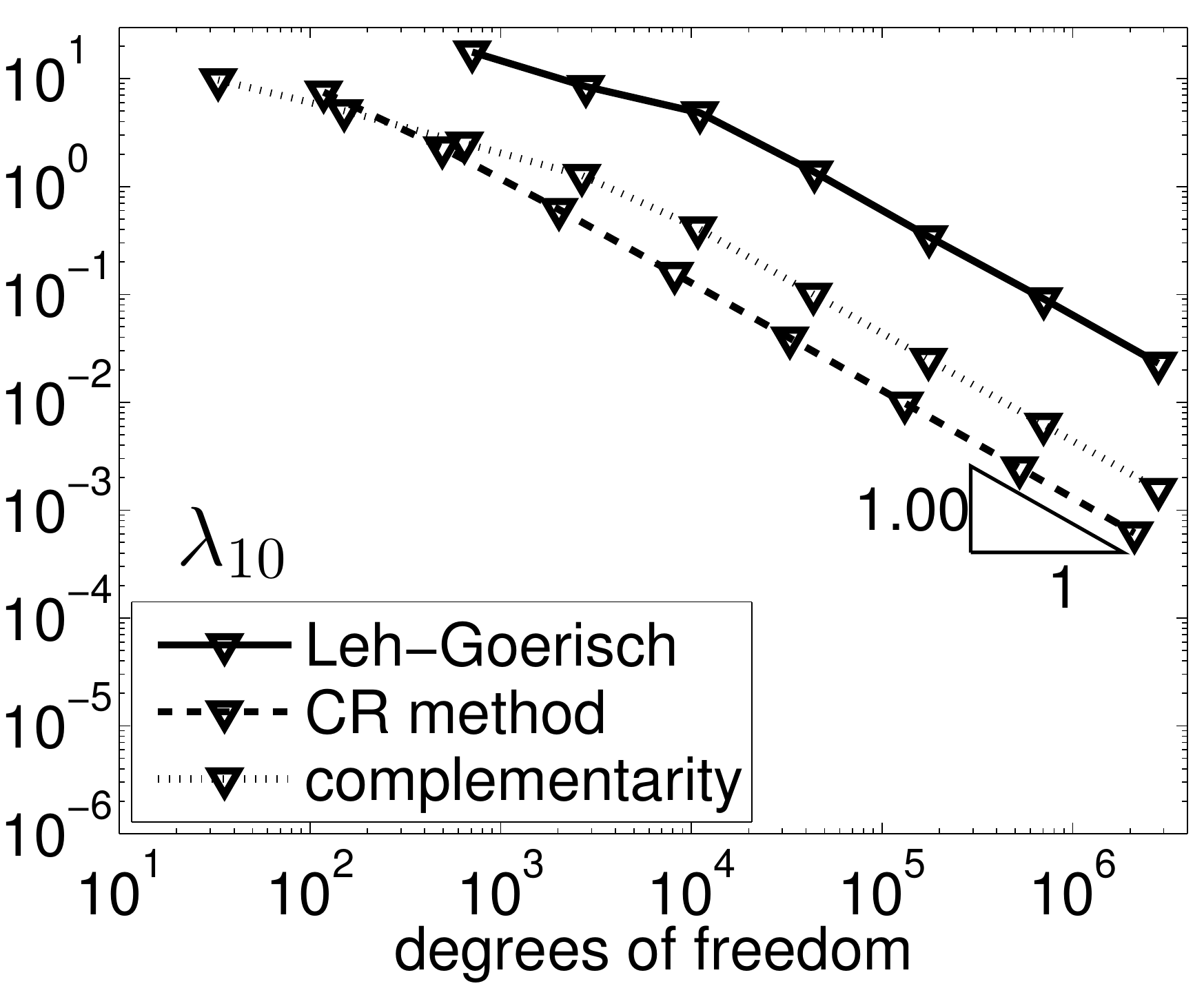}%
\caption{\label{fi:ctsqr_encl2}
Chopped off square.
Enclosure sizes of eigenvalues $\lambda_7$, $\lambda_8$, \dots, $\lambda_{10}$
obtained by  the lowest order Lehmann--Goerisch, CR, and complementarity methods. 
Triangles indicate experimental orders of convergence. 
}
\end{figure}

The three methods are compared using the same methodology as in the above sections. The initial mesh is illustrated in Figure~\ref{fi:sq_initmesh} (right). This time there is no obvious symmetry in the problem and the mesh is not symmetric.
For the Lehmann--Goerisch method we succeeded to refine the mesh 6 times, for the CR method 7 times, and for the complementarity method 8 times.
The eigenvalue bounds obtained on the finest available meshes are presented Table~\ref{ta:ctsqr_eig}.
Convergence curves of the lowest order versions of all three methods on a series of uniformly refined meshes are reported in Figures~\ref{fi:ctsqr_encl1}--\ref{fi:ctsqr_encl2}.

Resulting enclosure sizes are similar to those for the square, see Section~\ref{se:sqrpipi}. The complementarity and CR methods produce the smallest eigenvalue enclosures for the given number of degrees of freedom. The Lehmann--Goerisch method closely follows due to reasons mentioned above.

\section{Numerical results -- higher order}
\label{se:pdeg}

Both the Lehmann--Goerisch and complementarity methods have natural higher order versions.
This section compares their numerical performance for orders $k=1,2, \dots, 5$.
The CR method is not included, because its higher order version is not available.

\begin{figure}
\includegraphics[width=0.50\textwidth]{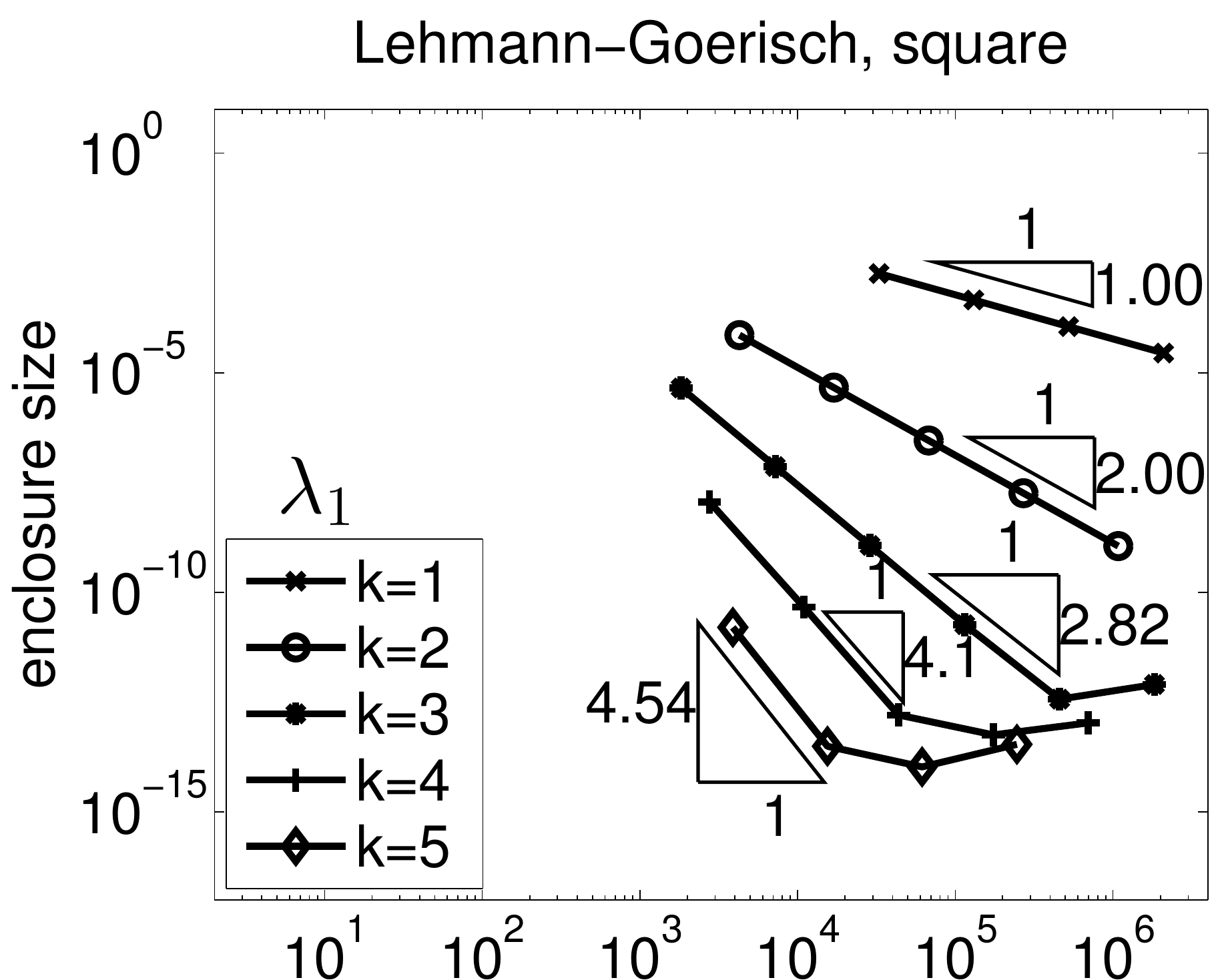}%
\includegraphics[width=0.48\textwidth]{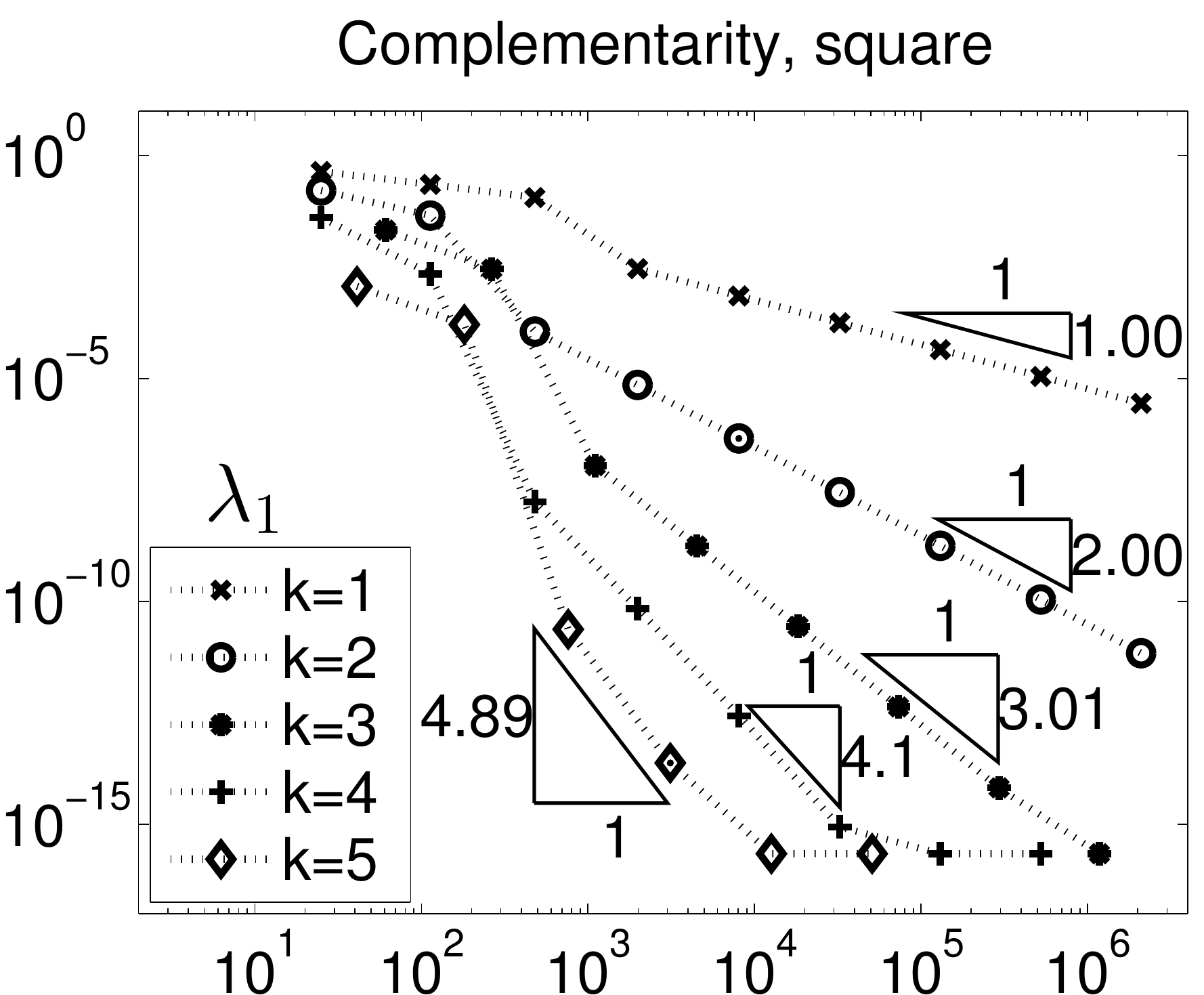}%
\\
\includegraphics[width=0.50\textwidth]{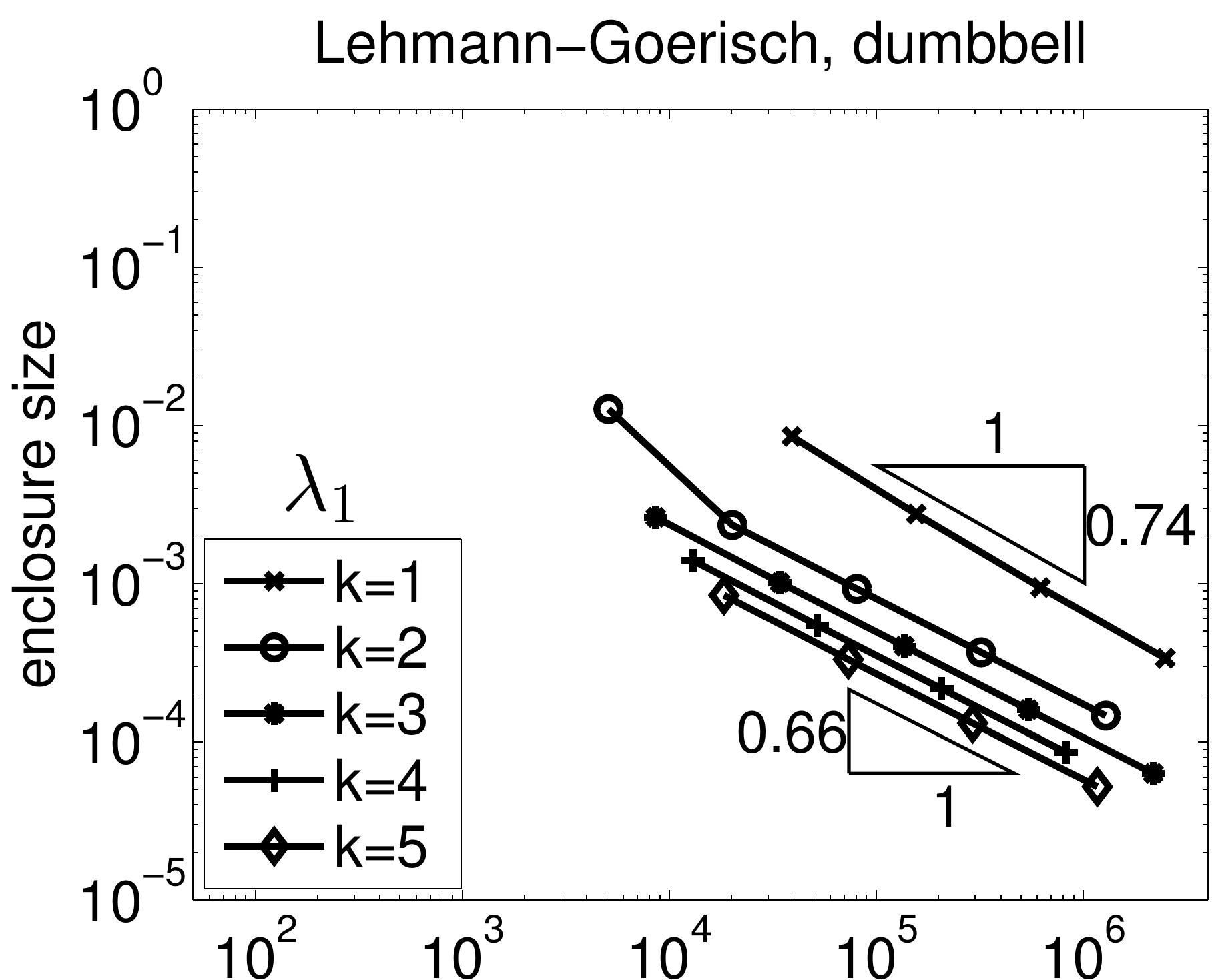}%
\includegraphics[width=0.48\textwidth]{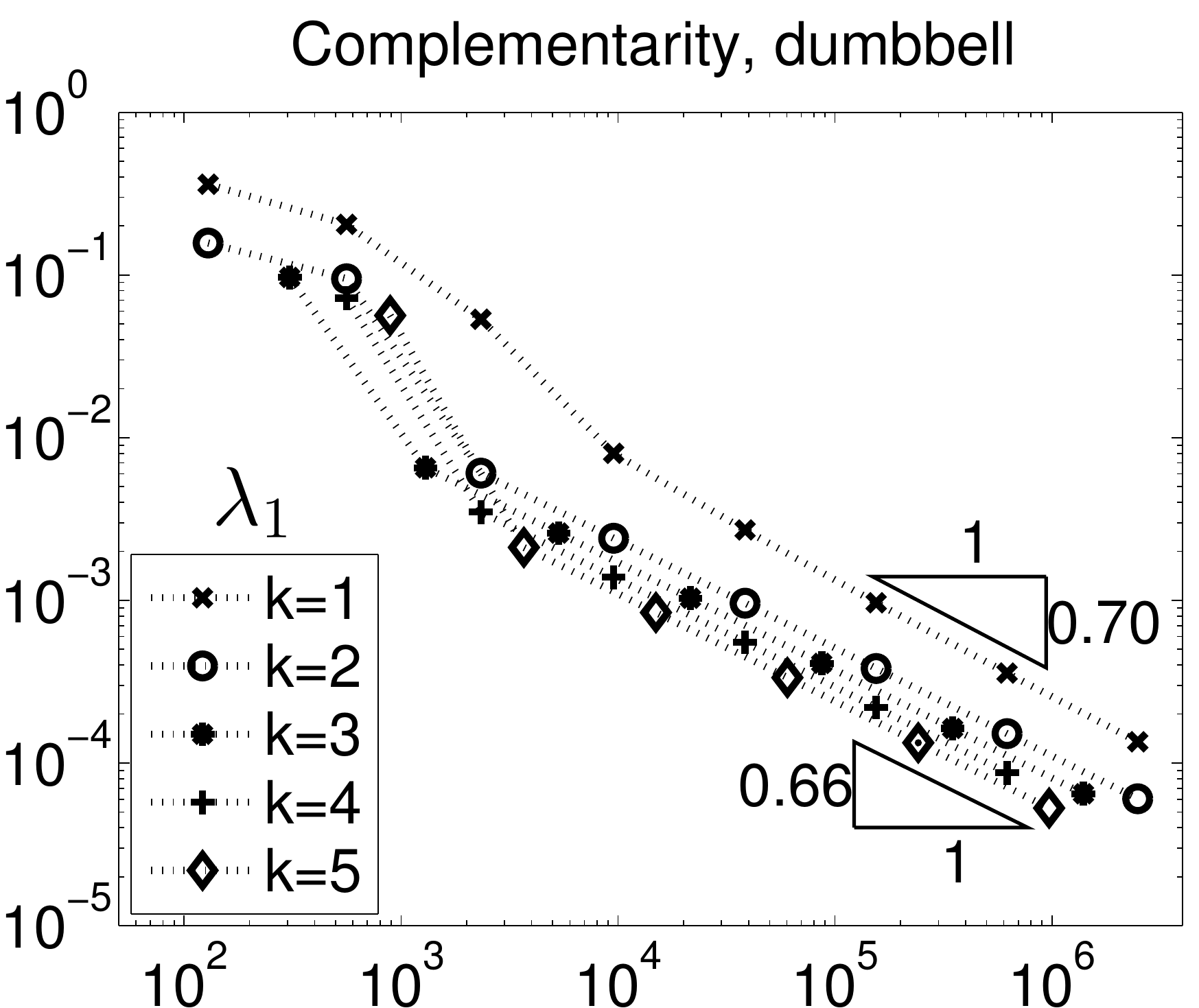}%
\\
\includegraphics[width=0.50\textwidth]{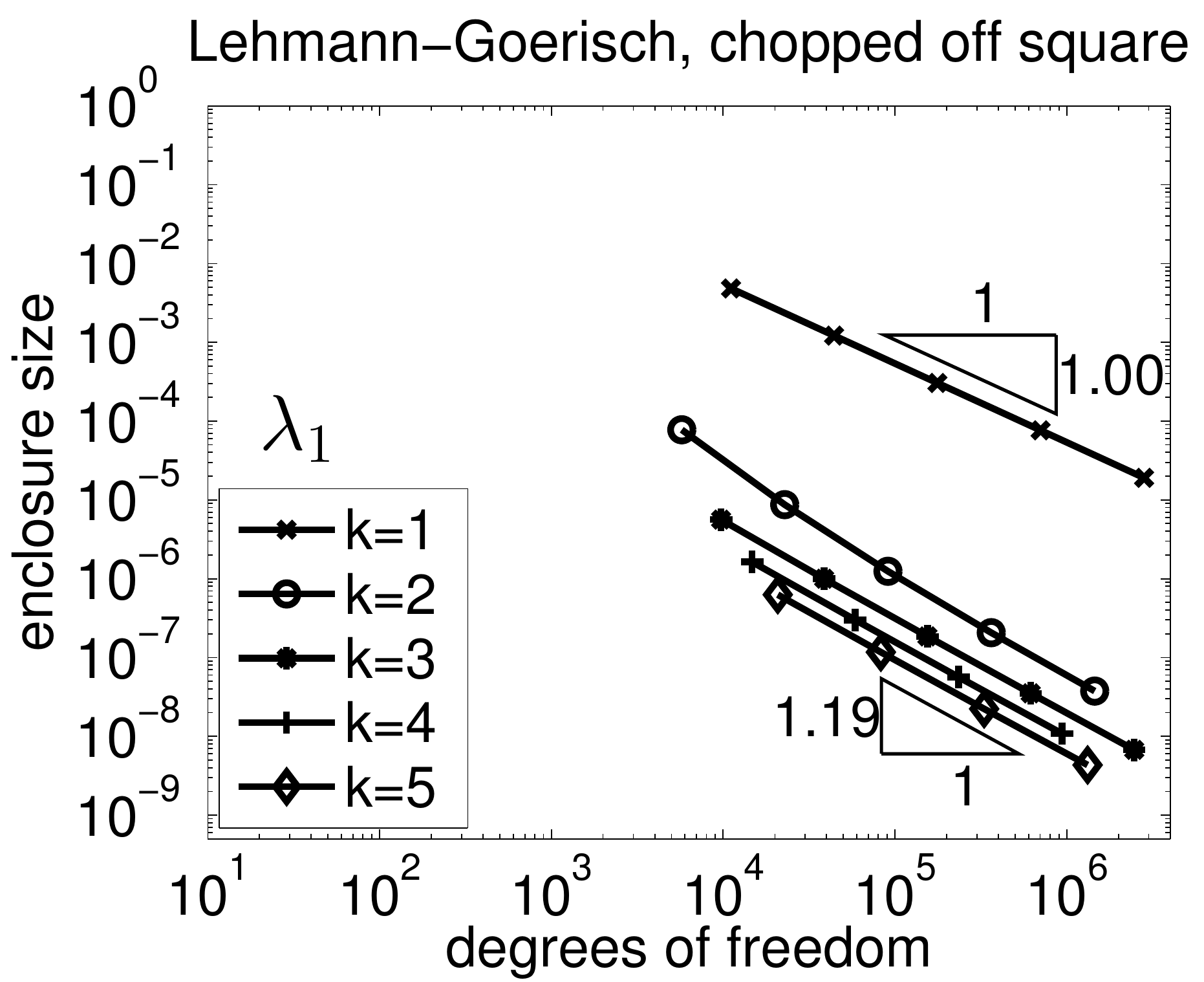}%
\includegraphics[width=0.48\textwidth]{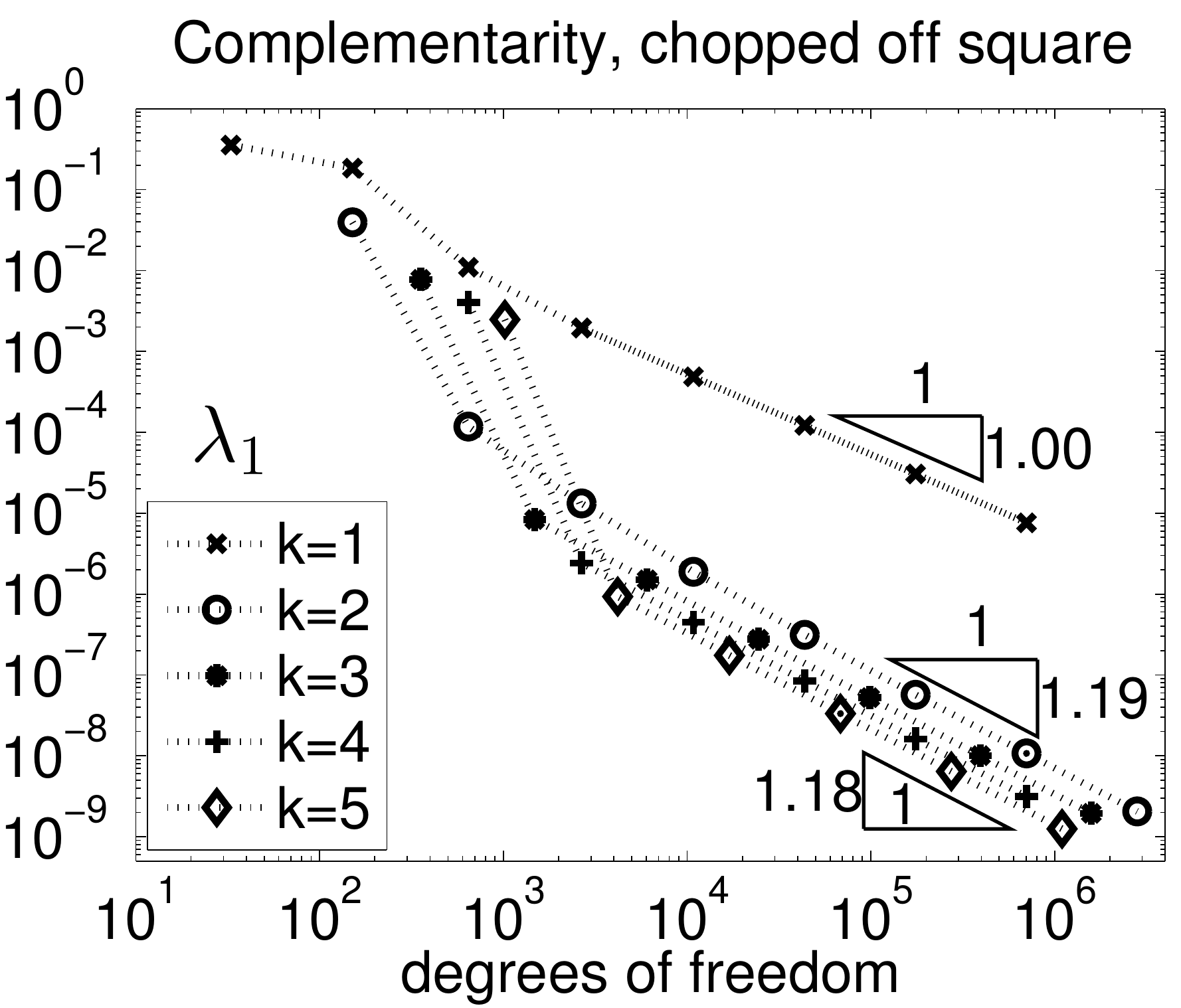}%
\label{fi:pdegl1}
\caption{Comparison of higher order Lehmann--Goerisch (left column) and complementarity (right column)
methods for the square (top row), dumbbell shaped domain (middle row), and chopped off square (bottom row).
Graphs show the dependence of enclosure sizes of the first eigenvalue with respect to the number of degrees of freedom.
}
\end{figure}

Figure~\ref{fi:pdegl1} compares convergence curves of the first eigenvalue with respect to the number of degrees of freedom for the square, dumbbell shaped domain, and the chopped off square. Domains, uniformly refined meshes, \emph{a priori} known lower bounds, and other parameters 
are the same as in Sections~\ref{se:sqrpipi}--\ref{se:ctsqr}, except for the order $k$.

The observed experimental orders of convergence for the square domain are close to the expected optimal rates for both methods. Enclosure sizes computed by the Lehmann--Goerisch method of order $k=3,4,5$ hit the limiting accuracy around $10^{-13}$, where inaccuracies of the algebraic solver of the large saddle point problem \eqref{eq:sig1}--\eqref{eq:sig2} probably start to dominate. Inaccuracies of the algebraic solver of local problems \eqref{eq:locprob1}--\eqref{eq:locprob2} are smaller and enclosure sizes computed by the complementarity method reach the machine precision for $k=3,4,5$.

For the dumbbell shaped domain the optimal rates of convergence are not achieved due to singularities of eigenfunctions in reentrant corners. Experimental rates of convergence of enclosure sizes computed by both methods are around $0.7$ for all $k=1,2,\dots,5$. In spite of this fact higher order versions produce more accurate results than lower order versions.

Concerning the square with chopped off corner, the optimal rate of convergence for both the Lehmann--Goerisch and complementarity methods is observed in the first order case only. Higher order versions of both methods exhibit experimental rates of convergence around $1.2$. This is caused by the loss of regularity of exact eigenfunctions in obtuse corners,
see \cite{Grisvard_prob_lim_polyg_86} where the regularity of the solution of Laplace equation is shown to decrease as the size of the largest angle of the polygonal domain increase. 

Optimal rates of convergence for general polygonal domains could be achieved by employing a mesh adaptive algorithm. All presented methods can be combined with an adaptive algorithm, but the corresponding numerical comparison is beyond the scope of this paper. 
We only note that localized version of quantity \eqref{eq:eta} computed in the complementarity method can be directly used to guide the adaptive mesh refinement \cite{VejSeb2017}.

Results for higher eigenvalues are similar. To illustrate the accuracy, Table~\ref{ta:pdeg} presents eigenvalue bounds computed by the fifth order Lehmann--Goerisch and complementarity methods on the finest meshes for the dumbbell shaped domain and the chopped off square. The corresponding bounds for the square are all computed to the level comparable to the machine precision and we do not present them.
We observe that both methods achieve comparable accuracy even for higher eigenvalues.


\begin{table}
\centerline{
\begin{tabular}{c|ccc|ccc}
 & \multicolumn{3}{c}{dumbbell shaped domain}& \multicolumn{3}{|c}{chopped off square} \\
 & lower bound   & lower bound & upper bound    & lower bound   & lower bound & upper bound    \\ 
 & Leh--Goerisch & complement. & Rayleigh--Ritz & Leh--Goerisch & complement. & Rayleigh--Ritz \\ \hline
$\lambda_{1}$ &   1.95576583 &   1.95575050 &   1.95580337 &  2.0042919800 &   2.0042919809 &   2.0042919821 \\
$\lambda_{2}$ &   1.96065834 &   1.96066662 &   1.96069147 &  5.0000350014 &   5.0000349937 &   5.0000350016 \\
$\lambda_{3}$ &   4.80050602 &   4.80044073 &   4.80080422 &  5.0301050037 &   5.0301050107 &   5.0301050229 \\
$\lambda_{4}$ &   4.82967840 &   4.82975419 &   4.82993162 &  8.0523670425 &   8.0523670504 &   8.0523670844 \\
$\lambda_{5}$ &   4.99682476 &   4.99682068 &   4.99683908 &  10.000502720 &   10.000502692 &   10.000502725 \\
$\lambda_{6}$ &   4.99683861 &   4.99684369 &   4.99685288 &  10.055330923 &   10.055330952 &   10.055330983 \\
$\lambda_{7}$ &   7.98680901 &   7.98678047 &   7.98697548 &  13.000742583 &   13.000742537 &   13.000742596 \\
$\lambda_{8}$ &   7.98687662 &   7.98694514 &   7.98704246 &  13.198058648 &   13.198058881 &   13.198058972 \\
$\lambda_{9}$ &   9.34821025 &   9.35022960 &   9.35732779 &  17.002654432 &   17.002654500 &   17.002654654 \\
$\lambda_{10}$&   9.50171274 &   9.50727405 &   9.51086516 &  17.064780396 &   17.064780837 &   17.064780947 \\
\end{tabular}
}
\caption{\label{ta:pdeg}
Eigenvalue bounds computed by the fifth order Lehmann--Goerisch and complementarity methods on the finest meshes for the dumbbell shaped domain and the chopped off square. The mesh sizes of the finest meshes 
$h_\mathrm{LG}$ and $h_\mathrm{cmpl}$ for the Lehmann--Goerisch and complementarity methods, respectively, were
$h_\mathrm{LG}=0.0736$ and $h_\mathrm{cmpl}=0.0368$ for the dumbbell shaped domain and $h_\mathrm{LG}=0.0224$ and $h_\mathrm{cmpl}=0.0112$ for the chopped off square.
}
\end{table}

\section{Conclusions}
\label{se:conclusions}

Two-sided bounds of the first ten eigenvalues of the Dirichlet Laplacian were computed for three numerical examples.
For this purpose, the Lehmann--Goerisch, CR, and complementarity methods were employed on a series of uniformly refined meshes and their accuracy was compared. 
All methods exhibit the optimal rate of convergence and the most accurate lower bounds with respect to the chosen numbers of degrees of freedom were obtained by the complementarity and CR methods for the lowest order case and by Lehmann--Goerisch and complementarity method in the higher order case.

For the lowest order case, the Lehmann--Goerisch method lags behind the other two methods in terms of the accuracy with respect to the number of degrees of freedom, because it requires to solve the large mixed finite element problem \eqref{eq:sig1}--\eqref{eq:sig2}. If fluxes $\bsigma_{h,i}$ were computed locally and efficiently as in the complementarity method then the convergence curves of the Lehmann--Goerisch method would be shifted considerably to the left, see Figures~\ref{fi:sq_encl1}--\ref{fi:ctsqr_encl2}, and the method would probably outperform the other two. A serious disadvantage of this method is the need of the \emph{a priori} lower bounds on eigenvalues. 
In \cite{BehMerPluWie2000} it is proposed to use the homotopy method to obtain them. The homotopy method, however, seems to be difficult to automatize in the context of the finite element method. On the other hand, an advantage of the Lehmann--Goerisch method is its generality and the fact that it is well developed in the literature. It can be generalized even to indefinite problems and to problems with essential spectrum. For example, its variant was used in \cite{Barrenechea2014} to find eigenvalue enclosures for the Maxwell eigenvalue problem. 
Another advantage of the Lehmann--Goerisch method is its straightforward generalization to higher order finite elements. If the computed eigenfunctions are smooth and do not contain singularities then the higher order approach provides high-precision results even on uniform meshes, see also \cite{LiuOis2013b}.

Concerning the CR method, we mention that
the dimension of the space $\VCR_h$ is larger than the dimension of $V_h$ with $k=1$, which is used as a base space for the lowest order version of the other two methods. Nevertheless, it is still considerably less than the number of degrees of freedom for the mixed finite element problem \eqref{eq:sig1}--\eqref{eq:sig2}.
The distinctive feature of the CR method is that it provides guaranteed lower bounds even on rough meshes and that
it does not require any \emph{a priori} information about the exact eigenvalues.
The \emph{a priori} information is not needed thanks to the explicitly known value of the interpolation constant $\kappa$. However, this constant is explicitly known only for specific operators, such as the Laplacian and biharmonic operator and for specific choices of finite elements, such as Crouzeix--Raviart and Morley elements.





The complementarity method is based on conforming finite elements and its lowest order version provides comparable accuracy as the CR method with respect to the number of degrees of freedom. Similarly as the Lehmann--Goerisch method, it requires an \emph{a priori} information about the exact eigenvalues.
% ...

Another similarity to the Lehmann--Goerisch method is that the complementarity method can be straightforwardly generalized to more complex problems. A difference is that flux reconstruction in the complementarity method is computed by solving small mixed finite element problems on patches of elements. These problems are independent and can be easily solved in parallel.

Higher order versions of both the Lehmann--Goerisch and complementarity methods provide highly accurate two-sided bounds on
eigenvalues. In the case of smooth eigenfunctions, the enclosure sizes computed by
the complementarity method reached the level of machine precision. Even in the
case of less regular eigenfunctions and uniformly refined meshes, the higher order
methods provide considerably higher accuracy than the lowest order methods.

To conclude, we are convinced that two-sided bounds of eigenvalues are highly relevant to compute, because they enable reliable control of the accuracy of the computed approximations. 
We believe that the presented methods are of practical value, because they are applicable in the context of the standard finite element method. We also believe that the presented results enable practitioners to choose the most suitable method for their purposes. Finally, we believe that these results encourage them to compute both upper and lower bounds on eigenvalues, because they enable full control of the accuracy and yield highly reliable numerical results.

\section*{Acknowledgements}
The author would like to thank anonymous referees for useful suggestions that lead to improvements of the paper, for proposing the example without obvious symmetries presented in Section~\ref{se:ctsqr}, and for proposing the comparison of higher order methods in Section~\ref{se:pdeg}.
Further, the author gratefully acknowledges the support of Neuron Fund for Support of Science, project no.~24/2016, and the institutional support RVO~67985840.

\bibliographystyle{amsplain} 
\bibliography{bibl}    

\end{document}